\documentclass[a4paper]{elsarticle}
\usepackage{lipsum}
\makeatletter
\def\ps@pprintTitle{%
 \let\@oddhead\@empty
 \let\@evenhead\@empty
 \def\@oddfoot{}%
 \let\@evenfoot\@oddfoot}
\makeatother
\usepackage[left=3cm,right=3cm,top=2cm,bottom=2cm,includeheadfoot]{geometry}
\usepackage{amsmath}
\usepackage{amsfonts}
\usepackage{amssymb}
\usepackage{amsfonts}
\usepackage{ MnSymbol}
\usepackage{amsmath}
\usepackage{cancel}
\usepackage{graphicx}
\usepackage[numbers]{natbib}

\usepackage[ansinew]{inputenc}
\usepackage{amsthm}
\usepackage[arrow, matrix, curve]{xy}
\usepackage{hyperref}
\newtheoremstyle{style}   
 {0.5cm}                 
 {0.5cm}                 
  {}                         
  {}                         
  {\normalfont\bfseries}  
  {\normalfont }{   }
  {}
\newtheorem{thm}{Theorem}
\newtheorem{satz}[thm]{Theorem}
\newtheorem{Proposition}[thm]{Proposition}

\newtheorem{Lemma}[thm]{Lemma}
\newtheorem{beispiel}[thm]{Example}
\newtheorem{kor}[thm]{Corollary}
\theoremstyle{remark}
\newtheorem{bemerkung}{Remark}
\theoremstyle{plain}
\newtheorem{definition}{Definition}
\newcommand{\ph}{\varphi}
\newcommand{\R}{\mathbb{R}}

\newcommand{\C}{\mathbb{C}}

\newcommand{\De}{\Delta_{p,q}}

\newcommand{\Pe}{\mathcal{P}}

\newcommand{\g}{\mathfrak{g}}

\allowdisplaybreaks[1]
\begin{document}

\begin{frontmatter}

\title{{\textbf{Conformal superalgebras via tractor calculus}}}


\author{Andree Lischewski}

\address{Department of Mathematics, Humboldt University, Rudower Chausse 25, 12489 Berlin, Germany}
\date{today}
\begin{abstract}
We use the manifestly conformally invariant description of a Lorentzian conformal structure in terms of a parabolic Cartan geometry in order to introduce a superalgebra structure on the space of twistor spinors and normal conformal vector fields formulated in purely algebraic terms on parallel sections in tractor bundles. Via a fixed metric in the conformal class, one reproduces a conformal superalgebra structure which has been considered in the literature before. The tractor approach, however, makes clear that the failure of this object to be a Lie superalgebra in certain cases is due to purely algebraic identities on the spinor module and to special properties of the conformal holonomy representation. Moreover, it naturally generalizes to higher signatures. This yields new formulas for constructing new twistor spinors and higher order normal conformal Killing forms out of existing ones, generalizing the well-known spinorial Lie derivative. Moreover, we derive restrictions on the possible dimension of the space of twistor spinors in any metric signature.
\end{abstract}

\begin{keyword}
twistor spinors \sep tractor calculus \sep conformal superalgebra
\MSC[2010] 53B30 \sep 53A30 \sep 53C27 \sep 15A66  

\end{keyword}
\ead{lischews@math.hu-berlin.de}

\end{frontmatter}
\tableofcontents
\section{Introduction}
The study of supersymmetric field theories in physics literature naturally leads to the notion of a (Lie) superalgebra, cf. \cite{med,med2,of12,of121,of122}, being defined as follows: Let $\g=\g_0 \oplus \g_1$ be a $\mathbb{Z}_2-$graded $\mathbb{K}-$vector space. For a homogeneous element $X \in \g$, we let $|X|:=i$ if $X \in \g_i$. $\g$ together with a bilinear map $[\cdot, \cdot ] : \g \times \g \rightarrow \g$ is called a $(\mathbb{K}-)$superalgebra if
\begin{enumerate}
\item  $[\cdot, \cdot ] : \g_i \times \g_j \rightarrow \g_{i+j}$,
\item For homogeneous elements $X,Y \in \g$ it holds that $[X,Y]=-(-1)^{|X||Y|}[Y,X]$.
\end{enumerate}
If moreover the Jacobi identity
\begin{align}
[X,[Y,Z]]=[[X,Y],Z]+(-1)^{|X||Y|}[Y,[X,Z]] \label{grj}
\end{align}
holds for all homogeneous elements, we call $\g$ a Lie superalgebra. Classification results for simple Lie superalgebras can be found in \cite{nahm}. It has been found in \cite{jgeo,suads,med,med2,med3,of121} that some superalgebras naturally appear geometrically. To this end, let $(M,g)$ be a smooth, oriented and time-oriented Lorentzian spin manifold with spinor bundle $S^g$ admitting distinguished spinor fields, e.g. parallel spinors, geometric Killing spinors or spinors being parallel wrt. a connection that depends on more bosonic data of the background in question. By a well-known squaring map, cf. \cite{lei,bl}, each spinor gives rise to a vector field and the spinor field equation translates into natural properties of the associated vector, i.e. being parallel or Killing, for instance. Moreover, vector fields act naturally on spinors by the spinorial Lie derivative as considered in \cite{kos,ha96}. In this way, one obtains a superalgebra naturally associated to $(M,g)$ whose even and odd part consist of distinguished vector- and spinor fields.
The algebraic structure of these infinitesimal symmetries also becomes important within the classification of background geometries on curved space which support some (rigid) supersymmetry, as has been initiated in physics literature in recent years, cf. \cite{pes,ym,CCKS1,CCKS2,CCKS3,fes1,cas}.\\
\newline
There is a conformal analogue of this superalgebra construction which has first been studied in \cite{ha96} and recently has been refined in \cite{raj,med,med2,med3}.
To this end note that besides the Dirac operator $D^g$ on a Lorentzian spin manifold, there is a complementary conformally covariant differential operator acting on spinors, called the Penrose-or twistor operator, and elements of its kernel are equivalently characterized as solutions of the twistor equation
\begin{align*}
\nabla^{S^g}_X \ph + \frac{1}{n} X \cdot D^g \ph = 0 \text{ for }X \in TM.
\end{align*}
There are many geometric classification results for manifolds admitting twistor spinors, cf. \cite{bl,leihabil,bfkg}, and recently they also appeared in Fefferman constructions in parabolic geometry, cf. \cite{hs,hs1}, and in the construction of conformal superalgebras in physics literature, cf. \cite{raj,med,med3,CCKS1,CCKS2,CCKS3}.\\ 
Under further orientability assumptions on $(M,g)$ every twistor spinor $\ph$ defines an associated Dirac current $V_{\ph} \in \mathfrak{X}(M)$ which turns out to be causal and conformal, see \cite{lei,bl}, and at least in the Lorentzian case their zero sets coincide, i.e. $Z_{\ph} = Z_{V_{\ph}}$. 
On the space $\mathfrak{X}^{nc}(M)\oplus \text{ker }P^g$ of \textit{normal} conformal vector fields and twistor spinors \cite{ha96,raj} introduce brackets by setting:

\begin{equation}\label{alr}
\begin{aligned} 
\left[V,W\right] &:=[V,W]_{\mathfrak{X}(M)},\\
[V,\ph]&:= V \circ \ph, \\ 
[\ph,V]&:= -V \circ \ph, \\
[\ph_1,\ph_2]&:=V_{\ph_1,\ph_2}, 
\end{aligned}
\end{equation}
where $V,W \in \mathfrak{X}^{nc}(M), \ph \in \text{ker }P^g$. $V \circ \ph$ is the spinorial Lie derivative (cf. \cite{kos}). It is proved in \cite{ha96,raj} that $\g:=\mathfrak{X}^{nc}(M)\oplus \text{ker }P^g$ together with these brackets is a superalgebra which is in general no Lie superalgebra. It has earlier been observed in \cite{ha96} that also the space ${\g}^{ec}:=\mathfrak{X}^{c}(M)\oplus \text{ker }P^g$ of conformal vector fields and twistor spinors equipped with the same brackets turns out to be a superalgebra which in general is no Lie superalgebra. We will discuss later why we choose only normal conformal vector fields in the even part.\\
\cite{med,med3} relates the superalgebra defined by (\ref{alr}) to the (local) classification of Lorentzian conformal structures admitting twistor spinors from \cite{leihabil}. One finds that (\ref{alr}) does not define a Lie superalgebra in case that there is a Fefferman metric in the conformal class or a Lorentzian Einstein Sasaki metric or a local splitting into a special Einstein product in the sense of \cite{al2}. In these cases the odd-odd-odd Jacobi identity fails to hold but the situation can be remedied by inclusion of a nontrivial R-symmetry in the construction of the algebra.\\
\newline
The mentioned constructions of conformal superalgebras involving twistor spinors all fix a metric in the conformal class. In contrast to this, our aim is the construction of a superalgebra canonically associated to a conformal spin structure by making use of conformal tractor calculus as developed in \cite{cs,baju,leihabil,feh}, for instance.
To this end, we use the well-known description of a pseudo-Riemannian conformal structure $(M,c=[g])$ of signature $(p,q)$,  where $n=p+q \geq 3$, as a parabolic Cartan geometry $(\mathcal{P}^1, \omega^{nc})$ of type $(G=O(p+1,q+1),P)$, where $P \subset G$ is the stabilizer of some isotropic ray, in the sense of \cite{sharp,cs,baju}. It leads to a well-defined algebraic conformal invariant, being the conformal holonomy group $Hol(M,c)$. As no canonical connection for $(M,c)$ can be defined on a reduction of the frame bundle of $M$, the Cartan geometry in question arises via a procedure called the first prolongation of a conformal structure, which naturally identifies $Hol(M,c) \cong Hol(\omega^{nc})$ with a (class of conjugated) subgroup of $O(p+1,q+1)$. Conformal holonomy groups turn out to be interesting objects in their own right, cf. \cite{arm,leihabil,baju,ln,lst,alt}, for instance. Normal conformal vector fields are in this language equivalently characterized as sections of the bundle $\mathcal{P}^1 \times_P \Lambda^2 \R^{p+1,q+1}$ that are parallel wrt. the connection induced by $\omega^{nc}$.
Furthermore, \cite{lei,baju,leihabil} shows that the twistor equation admits a conformally invariant reinterpretation in terms of conformal Cartan geometries. In fact, there is a naturally associated vector bundle $\mathcal{S}$ for a conformal spin manifold $(M,c)$ of signature $(p,q)$ with fibre $\Delta_{p+1,q+1}$, the spinor module in signature $(p+1,q+1)$. On $\mathcal{S}$, a natural lift of the conformal Cartan connection $\omega^{nc}$ induces a covariant derivative such that parallel sections of $\mathcal{S}$ correspond to twistor spinors via a fixed metric $g \in c$. In other words, $(M,c)$ admits a twistor spinor for one - and hence for all - $g \in c$ iff the lift of $Hol(M,c)$ to the spin group $Spin^+(p+1,q+1)$ which double covers $SO(p+1,q+1)$ stabilizes a nonzero spinor. Using these Cartan techniques has lead to a complete local classification of Lorentzian conformal structures admitting twistor spinors in \cite{leihabil}.\\
\newline
We present in this language a manifestly conformally invariant construction of a superalgebra $\g=\g_0 \oplus \g_1$, consisting of parallel tractor 2-forms and parallel tractor spinors,
\begin{equation*}
 \boxed{
\begin{aligned}
 &(M^{1,n-1},c) && & \rightarrow & \g= \g_0 \oplus \g_1 \text{(real) Superalgebra,}& \\
 &(M^{1,n-1},c) &\text{ with \textit{special} holonomy}& & \rightarrow & \g=\g_0 \oplus \g_1 \text{ \textit{Lie} Superalgebra.}&
\end{aligned}
}
\end{equation*}
defined on the level of tractors only. Here, the various brackets are given in purely algebraic terms by the obvious bracket on skew-symmetric endomorphisms, the natural Clifford-action of 2-forms on spinors and the squaring of spinors to 2-forms in signature $(2,n)$, cf. \cite{leihabil}. One should compare this to the construction of a (Lie) superalgebra for Riemannian manifolds admitting geometric Killing spinors via algebraic operations on the metric cone as done in \cite{suads}. We verify in section \ref{cts} that the so constructed superalgebra satisfies all Jacobi identities except the odd-odd-odd-one which has to be checked in a case-by-case analysis.\\
As we shall see in section \ref{dem}, this approach reproduces the superalgebra (\ref{alr}) from \cite{raj,ha96,med} when we fix a metric in the conformal class, which identifies parallel sections with conformal vector fields and twistor spinors, and thus it yields an equivalent description of the conformal symmetry superalgebra. However, we prove in Theorem \ref{hola} that the tractor approach as presented here has the advantage of giving purely algebraic conditions in terms of conformal holonomy exhibiting when the construction actually leads to a Lie superalgebra. \\
Furthermore, we present in section \ref{rsymme} the construction of a Lie superalgebra naturally associated to a Fefferman spin space via the inclusion of nontrivial R-symmetries on the tractor level. Again, the construction is purely algebraic and reproduces results of \cite{med} for a fixed metric in the conformal class.\\
\newline
Consequently, the tractor approach to conformal superalgebras induced by twistor spinors is manifestly conformally invariant, yields direct relations to conformal holonomy and shows that the failure of being a Lie superalgebra is due to purely algebraic identities on the level of spin tractors.\\
Furthermore, we see in section \ref{highersgn} that the tractor approach can also be used to generalize the whole construction to non-Lorentzian signatures. Doing this, one faces an immediate problem: In general, the map $\chi \mapsto \alpha^2_{\chi}$ mapping a spinor to the associated 2-form is nontrivial only in case $p+1=2$, i.e. Lorentzian signature. In arbitrary signature, a nontrivial map can be obtained by forming $\alpha^{p+1}_{\chi}$. However, there is no obvious natural generalization of the Lie bracket on $\Lambda^k_{p+1,q+1}$ for $k>2$. Nevertheless, we introduce a natural superalgebra structure on the space of parallel forms and twistor spinors formulated in a purely algebraic way, check Jacobi identities and describe the algebra wrt. a given metric in the conformal class. Surprisingly, one finds that in arbitrary signature only 2 of the 4 Jacobi identities need to be satisfied. This is illustrated by considering generic twistor spinors in signature $(3,2)$ (cf. \cite{hs1}) as an example. As a second example, we specialize the construction to special Killing forms and geometric Killing spinors on pseudo-Riemannian manifolds.
In all these cases one obtains by fixing a metric in the conformal class interesting new formulas in Propositions \ref{prr} and \ref{na} which produce new twistor spinors and conformal Killing forms out of existing ones and which can be viewed as generalizations of the spinorial Lie derivative.\\
Finally, we relate the dimension of the odd part $\g_1$ to special geometric structures in the conformal class:
In physics, one is often not only interested in the existence of solutions of certain spinor field equations, but wants to relate the existence of a certain number of maximally linearly independent solutions to local geometric structures, cf. \cite{of12,jhom,of122}. From a more mathematical perspective, \cite{cortes} studies the relation between the existence of a certain number of parallel-, Killing- and twistor spinors and underlying local geometries. We present conformal analogues of some of these results in section \ref{podi}. For instance, we show in Proposition \ref{stuv} that a pseudo-Riemannian manifold admitting more than $\frac{3}{4}$ of the maximal number of linearly independent twistor spinors is already conformally flat.\\
\newline
This article is organized as follows: In section \ref{srw} we provide the necessary ingredients from conformal spin geometry and its conformally invariant reformulation in terms of tractors. Section \ref{cts} introduces the conformal symmetry superalgebra in terms of tractors for Lorentzian manifolds and studies elementary properties whereas section \ref{dem} relates this construction to previous results from \cite{med,raj}. Section \ref{rsymme} elaborates on the construction of a tractor conformal superalgebra for Fefferman spaces via the inclusion of an R-symmetry whereas section \ref{66} applies the results obtained so far in low dimensions. In section \ref{highersgn} we leave the Lorentzian setting and show how the purely algebraic tractor-formulas generalize to arbitrary signatures. We conclude with some relations between the algebraic structure of the conformal symmetry algebra and local geometries in the conformal class in section \ref{podi}.

\section{Preliminaries from conformal spin geometry} \label{srw}

\subsection*{Relevant spinor algebra}
We consider $\R^{p,q}$, that is, $\R^n$, where $n=p+q$, equipped with the scalar product $\langle \cdot, \cdot \rangle_{p,q}$ of index $p$, given by $\langle e_i, e_j \rangle_{p,q} = \epsilon_i \delta_{ij}$, where $(e_1,...,e_n)$ denotes the standard basis of $\R^n$ and $\epsilon_{i \leq p} = -1 = -\epsilon_{i>p}$. Let $e_i^{\flat}:=\langle e_i, \cdot \rangle_{p,q} \in \left(\R^{p,q}\right)^*$. We denote by $Cl_{p,q}$ the Clifford algebra of $(\R^{n},- \langle \cdot, \cdot \rangle_{p,q})$ and by $Cl_{p,q}^{\C}$ its complexification. It is the associative real or complex algebra with unit multiplicatively generated by $(e_1,...,e_n)$ with the relations $e_ie_j+e_je_i=-2 \langle e_i,e_j \rangle_{p,q}$.\\
Let $Spin(p,q) \subset Cl(p,q)$ denote the spin group and $Spin^+(p,q)$ its identity component. There is a natural double covering $\lambda:Spin(p,q) \rightarrow SO(p,q)$ of the pseudo-orthogonal group. Restricting irreducible representations of $Cl(p,q)$ or $Cl^{\C}(p,q)$ (cf. \cite{har,lm}) leads to the real or complex spinor module $\Delta_{p,q}^{\R}$ resp. $\Delta_{p,q}^{\C}$, cf. \cite{ba81,har,lm}. Further,  $Cl_{p,q}^{(\C)}$ acts on $\De$ and as $\R^n \subset Cl_{p,q} \subset Cl^{\C}_{p,q}$, this defines the Clifford multiplication $\cdot$ of a vector by a spinor, which naturally extends to a multiplication by $k$-forms: Letting 
$\omega = \sum_{1 \leq i_1 <...< i_k \leq n} \omega_{i_1...i_k} e^{\flat}_{i_1} \wedge...\wedge e^{\flat}_{i_k} \in \Lambda^k_{p,q}:=\Lambda^k \left(\R^{p,q}\right)^*$ and $\ph \in \De$, we set
\begin{align} \omega \cdot \ph := \sum_{1 \leq i_1 <...< i_k \leq n} \omega_{i_1...i_k} e_{i_1} \cdot...\cdot  e_{i_k} \cdot \ph \in \De. \label{clform} \end{align}
$\Delta_{p,q}$ admits a $Spin^+(p,q)$ nondegenerate invariant inner product $\langle \cdot, \cdot \rangle_{\De}$ such that \begin{align} \langle X \cdot u, v \rangle_{\De} + (-1)^p \langle u, X \cdot v \rangle_{\De} = 0. \label{fg}\end{align} for all $u,v \in \De$ and $X \in \R^n$.
In the complex case, it is Hermitian, whereas in the real case it is symmetric if $p=0,1$ mod $4$ with neutral signature ($p \neq0$ and $q \neq 0$) or it is definite ($p=0$ or $q=0$). In case $p=2,3$ mod $4$, the pair $(\De^{\R},\langle \cdot , \cdot \rangle_{\De^{\R}})$ is a symplectic vector space.\\
\newline
There is an important decomposition of $\Delta_{p+1,q+1}$ into $Spin(p,q)-$modules. Let $(e_0,...,e_{n+1})$ denote the standard basis of $\R^{p+1,q+1}$. We introduce lightlike directions $e_{\pm} := \frac{1}{\sqrt{2}}(e_{n+1} \pm e_0)$.  One then has a decomposition 
\begin{align}
\R^{p+1,q+1} = \R e_- \oplus \R^{p,q}\oplus \R e_+ \label{sp1}
\end{align}
into $O(p,q)-$modules. We define the annihilation spaces $Ann(e_{\pm}):=\{ v \in \Delta_{p+1,q+1} \mid e_{\pm}\cdot v = 0 \}$. It follows that for every $v \in \Delta_{p+1,q+1}$ there is a unique $w \in \Delta_{p+1,q+1}$ such that $v=e_+ w + e_- w$, leading to a decomposition
\begin{align}
\Delta_{p+1,q+1} = Ann(e_+) \oplus Ann(e_-). \label{fs}
\end{align}
$Ann(e_{\pm})$ is acted on by $Spin(p,q) \hookrightarrow Spin(p+1,q+1)$ and there is an isomorphism $\chi: Ann(e_-) \rightarrow \Delta_{p,q}$ of $Spin(p,q)$-modules leading to the identification
\begin{equation}
\begin{aligned} \label{deco}
\Pi: {\Delta_{p+1,q+1}}_{|Spin(p,q)} & \rightarrow \Delta_{p,q} \oplus \Delta_{p,q}, \\
v=e_+w+e_-w & \mapsto (\chi(e_-e_+w),\chi(e_-w)) 
\end{aligned}
\end{equation}

Spinors are related to forms by squaring, cf. \cite{cor,nc}: For $n=r+s$ we\footnote{For the moment we change the notation from $(p,q)$ to $(r,s)$ because we will later apply these results in cases in conformal geometry, where $(r,s)=(p,q)$ and $(r,s)=(p+1,q+1)$.} define 
\begin{equation}\label{6}
\begin{aligned}
\Gamma^k: \Delta_{r,s} \times \Delta_{r,s} \rightarrow \Lambda^k_{r,s}\text{, }(\chi_1,\chi_2) \mapsto \alpha_{\chi_1,\chi_2}^k,\text{ where} \\
\langle \alpha_{\chi_1,\chi_2}^k,\alpha \rangle_{r,s} := d_{k,r} \left( \langle \alpha \cdot \chi_1, \chi_2 \rangle_{\Delta_{r,s}}\right) \textit{  } \forall \alpha \in \Lambda^k_{r,s}. 
\end{aligned}
\end{equation}
The map $d_{k,r}: \mathbb{K} \rightarrow \mathbb{K}$ is the identity for $\mathbb{K}=\R$, whereas for $\mathbb{K}=\C$ it is defined as follows: One finds for complex spinors $\chi \in \De^{\C}$ that $\langle \alpha \cdot \chi, \chi \rangle_{\Delta_{r,s}^{\C}}$ is either real or purely imaginary. This depends on $(r,s)$ and $k$ as well as the chosen representation and admissible scalar product, but not on $\chi$. One then chooses $d_{k,r} \in \{Re, Im\}$ so that $\alpha_{\chi}:=\alpha_{\chi,\chi}^k$is indeed a real form and -if possible- nontrivial. It is obvious that the algebraic Dirac form $\alpha_{\chi}^k$ is explicitly given by the formula
\begin{align}
\alpha_{\chi}^k =   \sum_{1 \leq i_1 < i_2 <...<i_l \leq n} \epsilon_{i_1}...\epsilon_{i_l} \cdot d_{k,r} \left(\langle e_{i_1}\cdot...e_{i_l}\cdot \chi, \chi  \rangle_{\Delta_{r,s}}\right) e^{\flat}_{i_1} \wedge...\wedge e^{\flat}_{i_l}. \label{bla}
\end{align}
For $k=1$ the vector $V_{\chi}:=\left(\alpha^1_{\chi}\right)^{\flat}$ is the Dirac current. The construction is nontrivial at least for $k=r$ since $\alpha^r_{\chi}=0 \Leftrightarrow \chi = 0$.

\subsection*{The twistor equation on spinors}
Let $(M,g)$ be a space- and time oriented, connected pseudo-Riemannian spin manifold of  index $p$ and dimension $n=p+q \geq 3$. By $\Pe^g$ we denote the $SO^+(p,q)$-principal bundle of all space-and time-oriented pseudo-orthonormal frames. A spin structure  of $(M,g)$ is then given by a $\lambda-$reduction $(\mathcal{Q}^g,f^g)$ of $\Pe^g$ to  $Spin^+(p,q)$. The associated bundle $S^g:=\mathcal{Q}^g \times_{Spin(p,q)} \De$ is called the real or complex spinor bundle. Its elements are classes $[u,v]$. Fibrewise application of spinor algebra defines Clifford multiplication $\mu: T^*M \otimes S^g \rightarrow S^g$ and the Levi Civita connection  on $(M,g)$ lifts via $df^g$ and $\lambda^*$ to a connection $\widetilde{\omega}^g \in \Omega^1(\mathcal{Q}^g,\mathfrak{spin}(p,q))$ which in turn induces a covariant derivative $\nabla^{S^g}$ on $S^g$, locally given by the formula 
\begin{align*}
\nabla^{S^g}_X \ph = X(\ph) + \frac{1}{2} \sum_{1 \leq k < l \leq n} \epsilon_i \epsilon_j g(\nabla^g_X s_k,s_l) s_ks_l \cdot \ph, 
\end{align*}
for $\ph \in \Gamma(S^g)$ and $X \in \mathfrak{X}(M)$, where $s=(s_1,...,s_n)$ is any local pseudo-orthonormal frame. The composition of $\nabla^{S^g}$ with Clifford multiplication defines the Dirac operator $D^g : \Gamma(S^g) \rightarrow \Gamma(S^g)$, whereas performing $\nabla^{S^g}$ followed by orthogonal projection onto the kernel of Clifford multiplication gives rise to the twistor operator 
\[ P^g : \Gamma(S^g)  \stackrel{\nabla^{S^g}}{\rightarrow} \Gamma(T^*M \otimes S^g ) \stackrel{g}{\cong} \Gamma(TM \otimes S^g)  \stackrel{\text{proj}_{\text{ker}\mu}}{\rightarrow} \Gamma(\text{ker } \mu). \]
Spinor fields $\ph \in \text{ker }P^g$ are called twistor spinors and they are equivalently characterized as solutions of the twistor equation
\[\nabla^{S^g}_X \ph + \frac{1}{n} X \cdot D^g \ph = 0 \text{    for all } X \in \mathfrak{X}(M). \]
$P^g$ is conformally covariant: Letting $\widetilde{g}=e^{2 \sigma}g$ be a conformal change of the metric, it holds (cf. \cite{bfkg}) that $P^{\widetilde{g}} \widetilde{\varphi} = e^{-\frac{\sigma}{2}} \left( P^g(e^{-\frac{\sigma}{2}} \varphi) \right) \widetilde{}$. In particular, $\varphi \in \Gamma(S^g)$ is a twistor spinor with respect to $g$ if and only if the rescaled spinor $e^\frac{\sigma}{2} \widetilde{\varphi} \in \Gamma(S^{\widetilde{g}})$ is a twistor spinor with respect to $\widetilde{g}$, where  $\widetilde{}:S^g \rightarrow S^{\widetilde{g}}$ denotes the natural identification of the spinor bundles, see \cite{ba81}. 

\subsection*{Conformally invariant formulation in terms of tractors}
As twistor spinors are in fact objects of conformal geometry, \cite{baju,leihabil,cs,feh} has developed a concept describing twistor spinors if one is only given a conformal class $c=[g]$ instead of a single metric $g \in c$. As a preparation for this, recall that for $G$ an arbitrary Lie group with closed subgroup $P$ a Cartan geometry of type $(G,P)$ on a smooth manifold $M$ of dimension dim$(G/P)$ is specified by the data $(\mathcal{G} \rightarrow M, \omega)$, where $\mathcal{G}$ is a $P-$principal bundle over $M$ and $\omega \in \Omega^1(\mathcal{G},\mathfrak{p})$, called the Cartan connection, is $Ad$-equivariant wrt. the $P-$action, reproduces the generators of fundamental vector fields and gives a pointwise linear isomorphism $T_u \mathcal{G} \cong \mathfrak{g}$. The $P-$bundle $G \rightarrow G/P$ together with the Maurer-Cartan form of $G$ serves as flat and homogeneous model. As a Cartan connection does not allow one to distinguish a connection in the sense of a right-invariant horizontal distribution in $\mathcal{G}$, it is convenient to pass to the enlarged principal $G-$bundle $\overline{\mathcal{G}}:= \mathcal{G} \times_P G$ on which $\omega$ induces a principal bundle connection $\overline{\omega}$, uniquely determined by $\iota^* \overline{\omega} = \omega$, where $\iota: \mathcal{G} \hookrightarrow \overline{\mathcal{G}}$ is the canonical inclusion. For detailed introduction to Cartan geometries, we refer to \cite{sharp,cs}. \\
Applied to our setting, let $(M,c)$ be a connected, space- and time oriented conformal manifold of signature $(p,q)$ and dimension $n=p+q \geq 3$. It is well known that $c$ is equivalently, in the sense of \cite{cs}, encoded in a Cartan geometry $(\mathcal{P}^1 \rightarrow M, \omega^{nc})$ naturally associated to it via a construction called the first Prolongation of a conformal structure, cf. \cite{baju,cs}. In this case, the group $G$ is given by $G=SO^+(p+1,q+1)$ and the parabolic subgroup $P=Stab_{\R^+e_-}G$ is realized as the stabilizer of the lightlike ray $\R^+e_-$ under the natural $G-$action on $\R^{p+1,q+1}$. The homogeneous model is then given by $G/P \cong S^p \times S^q$ equipped with the obvious signature $(p,q)-$conformal structure. $\omega^{nc} \in \Omega^1(\mathcal{P}^1,\mathfrak{g})$ is called the normal conformal Cartan connection, and given $\mathcal{P}^1$, it s uniquely determined by the normalization condition $\partial^* \Omega^{nc}=0$ on its curvature $\Omega^{nc}:\mathcal{P}^1 \rightarrow Hom(\Lambda^2 \R^n, \mathfrak{so}(p+1,q+1)$, where $\partial^*$ denotes the Kostant codifferential, cf. \cite{cs}.\\
Given the standard $G-$action on $\R^{p+1,q+1}$, we obtain the associated standard tractor bundle $\mathcal{T}(M):=\mathcal{P}^1 \times_P \R^{p+1,q+1} = \overline{\mathcal{P}}^1 \times_G \R^{p+1,q+1}$ on which $\overline{\omega}^{nc}$ induces a covariant derivative $\nabla^{nc}$ which is metric wrt. the bundle metric $\langle \cdot, \cdot \rangle_{\mathcal{T}}$ on $\mathcal{T}(M)$ induced by the standard inner product on $\R^{p+1,q+1}$, and $\nabla^{nc}$ is therefore viewed as the conformal analogue of the Levi-Civita connection, making it reasonable to define the conformal holonomy of $(M,c)$ for $x \in M$ to be
\[ Hol_x(M,c):=Hol_x(\nabla^{nc}) \subset SO^+(\mathcal{T}_x(M), \langle \cdot , \cdot \rangle_{\mathcal{T}}) \cong SO^+(p+1,q+1). \]
By means of a metric in the conformal class, the conformally invariant objects introduced so far admit a more concrete description. Concretely, any fixed $g \in c$ induces a so-called Weyl-structure in the sense of \cite{cs} and leads to a $SO^+(p,q) \hookrightarrow G$-reduction $\sigma^g: \mathcal{P}^g \rightarrow \mathcal{P}^1$. Here, $\mathcal{P}^g$ denotes the orthonormal frame bundle for $(M,g)$. It follows with the decomposition (\ref{sp1}) that there is a $g-$metric splitting of the tractor bundle
\begin{align} 
\mathcal{T}(M) \stackrel{\Phi^g}{\cong} \underline{\R} \oplus TM \oplus \underline{\R}, \label{phig}
\end{align}
under which tractors correspond to elements $(\alpha,X,\beta)$ and the tractor metric takes the form
 \begin{align} \langle (\alpha_1, Y_1, \beta_1), (\alpha_2, Y_2, \beta_2) \rangle_{\mathcal{T}} = \alpha_1 \beta_2 + \alpha_2 \beta_1 + g(Y_1,Y_2). \label{bum} \end{align}
The metric description of the tractor connection $\nabla^{nc}$, i.e. $\Phi^g \circ \nabla^{nc} \circ (\Phi^g)^{-1}$ is (cf. \cite{baju})
 \begin{align} \nabla_X^{nc} \begin{pmatrix} \alpha \\ Y \\ \beta \end{pmatrix} = \begin{pmatrix} X(\alpha) + K^g(X,Y) \\ \nabla_X^g Y + \alpha X - \beta K^g(X)^{\sharp} \\ X(\beta) - g(X,Y) \end{pmatrix}, \label{trad} \end{align}
where $K^g := \frac{1}{n-2} \cdot \left( \frac{scal^g}{2(n-1)}  \cdot g - Ric^g  \right)$ is the Schouten tensor.\\
\newline 
Conformal Cartan geometry allows a conformally invariant construction of the twistor operator $P^g$. To this end, suppose that $(M,c)$ is additionally spin for one - and hence for all - $g \in c$.  Then the above construction admits a lift to a conformal spin Cartan geometry $(\mathcal{Q}^1, \widetilde{\omega}^{nc})$ of type $(\widetilde{G}:=Spin^+(p+1,q+1),\widetilde{P}:=\lambda^{-1}(P))$ with associated spin tractor bundle
\[ \mathcal{S}=\mathcal{S}(M):= \mathcal{Q}^1 \times_{\widetilde{P}} \Delta_{p+1,q+1}^{\R}, \]
on which $\mathcal{T}(M)$ acts by fibrewise Clifford multiplication and $\widetilde{\omega}^{nc}$ induces a covariant derivative $\nabla^{\mathcal{S}}$ on $\mathcal{S}$.
Fixing a metric $g \in c$ leads to a $Spin^+(p,q) \hookrightarrow Spin^+(p+1,q+1)$-reduction $\widetilde{\sigma}^g: \mathcal{Q}^g \rightarrow \mathcal{Q}^1$ which covers $\sigma^g$. We let $\overline{\mathcal{Q}^1}$ denote the enlarged $Spin^+(p+1,q+1)$-principal bundle and use $g$ to identify $\mathcal{S}(M) \cong Q^g_+ \times_{Spin^+(p,q)} \Delta_{p+1,q+1}$.
Together with the isomorphism (\ref{deco}), this leads to the $g-$metric identification
\begin{equation} \label{gdg}
\begin{aligned}
\widetilde{\Phi}^g: \mathcal{S}(M) &\rightarrow S^g(M) \oplus S^g(M), \\ 
[\widetilde{\sigma}^g(\widetilde{s}^g),v] &\mapsto [\widetilde{s}^g,\Pi(v)]
\end{aligned}
\end{equation}
with projections $proj^g_{\pm}$ to the annihilation spaces. One calculates that under (\ref{gdg}), $\nabla^{nc}$ is given by the expression (cf. \cite{baju})
\begin{align*}
\nabla^{nc}_X \begin{pmatrix} \ph  \\ \phi \end{pmatrix} = \begin{pmatrix} \nabla_X^{S^g} & -X \cdot \\ \frac{1}{2}K^g(X) \cdot & \nabla^{S^g}_X \end{pmatrix} \begin{pmatrix} \ph  \\ \phi \end{pmatrix}.
\end{align*}
As every twistor spinor $\ph \in \text{ker }P^g$ satisfies $\nabla^{S^g}_X \ph = \frac{n}{2}K(X) \cdot \ph$, cf. \cite{bfkg}, this yields a reinterpretation of twistor spinors in terms of conformal Cartan geometry. 
Namely for any $g \in c$, the vector spaces ker $P^g$ and parallel sections in $\mathcal{S}(M)$ wrt. $\nabla^{nc}$ are naturally isomorphic via
\begin{align*}
 \text{ker }P^g \rightarrow \Gamma(S^g(M) \oplus S^g(M)) \stackrel{\left(\widetilde{\Phi}^g\right)^{-1}}{\cong } \Gamma(\mathcal{S}(M))\text{,   } 
 \ph \mapsto \begin{pmatrix} \ph \\  -\frac{1}{n}D^g \ph \end{pmatrix} \stackrel{\left(\widetilde{\Phi}^g\right)^{-1}}{\mapsto} \psi \in Par(\mathcal{S}_{\mathcal{T}}(M), \nabla^{nc}), 
 \end{align*}
 i.e. a spin tractor $\psi \in \Gamma(\mathcal{S}(M))$ is parallel iff for one (and hence for all $g \in c$), it holds that $\ph:= \widetilde{\Phi}^g({proj}_+^g \psi) \in \text{ker }P^g$ and in this case $D^g \ph = -n \cdot \widetilde{\Phi}^g({proj}_-^g \psi)$.\\
In terms of conformal holonomy, the space of twistor spinors is thus in bijective correspondence to the space of spinors fixed by the lift of the conformal holonomy representation to $Spin^+(p+1,q+1)$, i.e. in the simply-connected case we have for $x \in M$ that
\begin{align} \label{sytr}
\text{ker }P^g \cong \{v \in \mathcal{S}_x \cong \Delta_{p+1,q+1} \mid \lambda_*^{-1}(\mathfrak{hol}_x(M,[g])) \cdot v = 0 \}.
\end{align}

\subsection*{The twistor equation on forms}
There is a canonical way of associating other parallel tractors to a twistor spinor. To this end, we introduce the tractor $(k+1)$-form bundle $\Lambda^{k+1}_{\mathcal{T}}(M):= \mathcal{P}^1 \times_P \Lambda^{k+1}_{p+1,q+1}$ on which again $\omega^{nc}$ induces a covariant derivative $\nabla^{nc} : \Gamma(\Lambda^{k+1}_{\mathcal{T}} (M)) \rightarrow \Gamma(T^{*}M \otimes \Lambda^{k+1}_{\mathcal{T}} (M))$. Fixed $g \in c$ allows us to describe tractor forms in terms of usual differential forms with the help of the following algebraic construction, using the decomposition $\R^{p+1,q+1} \cong \R e_- \oplus \R^{p,q} \oplus \R e_+$. Clearly, every form $\alpha \in \Lambda_{p+1,q+1}^{k+1}$ decomposes into
\begin{align}
\alpha = e_+^{\flat} \wedge \alpha_+ + \alpha_0 + e^{\flat}_- \wedge e_+^{\flat} \wedge \alpha_{\mp} + e_-^{\flat} \wedge \alpha_- \label{mind}
\end{align}
for uniquely determined forms $\alpha_-,\alpha_+ \in \Lambda^k_{p,q}, \alpha_0 \in \Lambda^{k+1}_{p,q}$ and $\alpha_{\mp } \in \Lambda^{k-1}_{p,q}$. Using this decomposition, the restriction of the standard action $O(p+1,q+1) \rightarrow GL\left(\Lambda^{k+1}_{p+1,q+1}\right)$ to $O(p,q) {\hookrightarrow} O(p+1,q+1)$  defines an isomorphism of $O(p,q)$-modules,
\begin{align*} \Lambda^{k+1}_{p+1,q+1} \cong \Lambda^k_{p,q} \oplus \Lambda^{k+1}_{p,q} \oplus \Lambda^{k-1}_{p,q} \oplus \Lambda^k_{p,q}.
\end{align*} 

This gives the $g$-metric representation of the tractor $(k+1)$-form bundle:
\begin{align*}
\Phi_{\Lambda}^g: \Lambda^{k+1}_{\mathcal{T}} (M) \stackrel{g}{\rightarrow} \Lambda^{k}(M) \oplus \Lambda^{k+1}(M) \oplus \Lambda^{k-1}(M) \oplus \Lambda^{k}(M).
\end{align*}
Applying this pointwise yields that each tractor $(k+1)$-form $\alpha \in \Omega^{k+1}_{\mathcal{T}}(M):=\Gamma\left(\Lambda^{k+1}_{\mathcal{T}}(M)\right)$ uniquely corresponds via $g \in c$ to a set of differential forms, 
\begin{align} \Phi_{\Lambda}^g\left( \alpha \right) = (\alpha_+ , \alpha_0 , \alpha_{\mp}, \alpha_-) \in  \Omega^k(M) \oplus \Omega^{k+1}(M) \oplus  \Omega^{k-1}(M) \oplus  \Omega^{k}(M). \label{dr} 
\end{align}
We further introduce the $g$-dependent projections
\begin{align*}
proj^g_{\Lambda,+}: \Omega_{\mathcal{T}}^{k+1}(M) &\rightarrow \Omega^k(M) \\
\alpha & \mapsto \alpha_+\text{, where } \Phi_{\Lambda}^g\left( \alpha \right) = (\alpha_+ , \alpha_0 , \alpha_{\mp}, \alpha_-)
\end{align*}
The operator $ \Phi_{\Lambda}^g \circ \nabla^{nc} \circ \left( \Phi_{\Lambda}^g \right)^{-1}$ satisfies
\begin{align} \label{stu}
\nabla_X^{nc} \alpha \stackrel{g}{=} \begin{pmatrix} \nabla^{g}_X & - X \invneg & -X^{\flat} \wedge & 0 \\ -K^g(X) \wedge &  \nabla^{g}_X & 0 & X^{\flat} \wedge \\ - \left(K^g(X)\right)^{\sharp} \invneg & 0 & \nabla^{g}_X & -X \invneg \\ 0 & \left(K^g(X)\right)^{\sharp} \invneg & -K^g(X) \wedge & \nabla^{g}_X \end{pmatrix} \begin{pmatrix} \alpha_+ \\ \alpha_0 \\ \alpha_{\mp} \\ \alpha_- \end{pmatrix}.
\end{align}
Finally, let $\alpha \in \Omega^{k+1}_{\mathcal{T}}(M)$ be a tractor $(k+1)-$form on $(M,c)$. Fix $g \in c$ and $\widetilde{g}= e^{2 \sigma} g \in c$ and let $\alpha_+={proj}^g_{\Lambda,+} \alpha$, $\widetilde{\alpha}_+={proj}^{\widetilde{g}}_{\Lambda,+} \alpha \in \Omega^k(M)$. These forms are related by (cf. \cite{nc}) 
\begin{align}
\widetilde{ \alpha}_+ = e^{(k+1) \sigma} \alpha_+. \label{wit}
\end{align}

The link between the above tractor forms and twistor spinors is given as follows: First, let $\ph_{1,2} \in \Gamma(S^g)$ and $\psi_{1,2} \in \Gamma(\mathcal{S}(M))$ be arbitrary spinor fields.  The algebraic construction (\ref{6}) can be made global by defining the following forms $\alpha_{\psi_1,\psi_2}^k \in \Omega^{k}_{\mathcal{T}}(M)$ and $\alpha^k_{\ph_1,\ph_2} \in \Omega^k(M)$ for every $k \in \mathbb{N}$:
\begin{equation} \label{3245}
\begin{aligned}
\begin{array}{llllllll}
\langle \alpha_{\psi_1,\psi_2}^k , \alpha \rangle_{\mathcal{T}} &:=& d_{k,p+1}& \left(\langle \alpha \cdot \psi_1 , \psi_2 \rangle_{\mathcal{S}} \right)& \forall & \alpha & \in & \Omega_{\mathcal{T}}^k(M), \\
g ( \alpha_{\ph_1,\ph_2}^k , \alpha ) &:=& d_{k,p}& \left(\langle \alpha \cdot \ph_1 , \ph_2 \rangle_{S^g} \right)& \forall & \alpha & \in & \Omega^k(M).
\end{array}
\end{aligned}
\end{equation}
It is straightforward to check that for $\psi_{1,2} \in Par (\mathcal{S}(M), \nabla^{nc})$ and $k \in \mathbb{N}$, the tractor $k-$form $\alpha^{k}_{\psi_1,\psi_2}$ is parallel wrt. $\nabla^{nc}$.\\
\newline
More generally, we call every parallel tractor $(k+1)$-form $\alpha \in Par(\Lambda^{k+1}_{\mathcal{T}}(M),\nabla^{nc}) \subset \Omega^{k+1}_{\mathcal{T}}(M)$, i.e $\nabla^{nc} \alpha = 0$, a twistor-$(k+1)$-form.
Using (\ref{stu}), \cite{nc} calculates that $\nabla^{nc} \alpha = 0$ implies for one - and hence for all - $g \in c$ the conformally covariant condition  
\begin{align} \label{soga}
\Phi^g_{\Lambda}(\alpha) = \left(\alpha_+,\frac{1}{k+1}d\alpha_+,-\frac{1}{n-k+1}d^*\alpha_+,\Box_k \alpha_+ \right),
\end{align}
whereby we have set 
\begin{align*} \Box_k :=\begin{cases} \frac{1}{n-2k} \left( - \frac{scal^g}{2(n-1)} + \nabla^* \nabla  \right), & n \neq 2k, \\ 
\frac{1}{n} \left(\frac{1}{k+1} (d^*d + dd^*)+\sum_{i=1}^n \epsilon_i \left(s_i \invneg (K^g(s_i)^{\flat} \wedge \cdot ) - s_i^{\flat} \wedge (K^g(s_i) \invneg \cdot )\right)\right), & n = 2k. \end{cases}
\end{align*}
Here, $s=(s_1,...,s_n)$ is a local section of $\Pe^g$ and $\nabla^*$ denotes the formal adjoint of $\nabla = \nabla^g$. \\
That is, $\nabla^{nc}\alpha = 0$ translates via (\ref{soga}) into a differential system for $\alpha_+$ only, and we call $\alpha_+ \in \Omega^k (M)$ arising in this way a normal conformal Killing $k$-form (or shorty, a nc-Killing form). We denote the set of these forms by $\Omega_{nc,g}^k(M)$ . Only considering the first equation in (\ref{stu}) leads to conformal Killing forms (cf. \cite{sem}). A conformal Killing form which is closed for some metric $g \in c$ is called a Killing form for $(M,g)$. In summary, each $g \in c$ leads to a natural isomorphism
\begin{align*}
Par\left(\Lambda^{k+1}_{\mathcal{T}}(M),\nabla^{nc}\right) \ni \alpha \mapsto {proj}_{\Lambda,+}^g (\alpha) \in \Omega^k_{nc,g}(M), \\
\end{align*}
where the inverse is given by $\alpha_+ \mapsto \left(\Phi^g_{\Lambda}\right)^{-1} \left(\alpha_+,\frac{1}{k+1}d\alpha_+,-\frac{1}{n-k+1}d^*\alpha_+,\Box_k \alpha_+ \right)$.\\
\newline
We turn again to twistor spinors. Let $\psi \in Par\left(\Lambda_{\mathcal{T}}^{k+1}(M)\right)$, $g \in c$ and $\ph:=\widetilde{\Phi}^g \left( {proj}^g_+ \psi \right) \in \text{ker }P^g$. It has been shown in \cite{ns} that there are constants $c^i_{k,p} \neq 0$ for $i=1,2$ such that
\begin{align}
{proj}_{\Lambda,+}^g \left(\alpha^{k+1}_{\psi} \right) = c^1_{k,p} \cdot \alpha^k_{\ph} \text{ and }{proj}_{\Lambda,-}^g \left(\alpha^{k+1}_{\psi} \right)=c^2_{k,p} \cdot \alpha^k_{D^g\ph}. \label{tuffi}
\end{align}
In particular, (\ref{tuffi}) reveals that for every twistor spinor $\ph \in \text{ker }P^g$, the forms $\alpha^k_{\ph}$ are nc-Killing forms. Together with the conformal transformation behaviour of $\ph$ and $\alpha_+$ under a change $\widetilde{g}=e^{2 \sigma}g$, this may be visualized in the following commutative diagram:
\begin{center}
\boxed{
\begin{aligned}
\begin{xy}
  \xymatrix{
      \ph \in \text{ker }P^g \ar[r]^{\left(\widetilde{\Phi}^g \circ {proj}^g_+ \right)^{-1}} \ar[dd]^{\text{nc-Killing}}    &   \psi \in Par(\mathcal{S}(M),\nabla^{nc}) \ar[r]^{\widetilde{\Phi}^{\widetilde{g}} \circ {proj}^{\widetilde{g}}_+} \ar[dd]^{\text{twistor form}} & e^{\frac{\sigma}{2}}\widetilde{\ph} \in \text{ker }P^{\widetilde{g}} \ar[dd]^{\text{nc-Killing}}\\
      & &  \\
      c^1_{k,p} \cdot\alpha^k_{\ph} \in \Omega_{nc,g}^k(M) \ar[r]^{ \left({proj}^g_{\Lambda,+} \right)^{-1}  }             &  \alpha^{k+1}_{\psi} \in \Omega^{k+1}_{\mathcal{T}}(M)\ar[r]^{{proj}^{\widetilde{g}}_{\Lambda,+}}  & c^1_{k,p}\cdot e^{(k+1)\sigma}\alpha_{\ph}^k \in \Omega_{nc,\widetilde{g}}^k(M)  \\
  }
\end{xy}
\end{aligned}}
\end{center}

For the special case of a twistor 2-form $\alpha \in \Omega^2_{\mathcal{T}}(M)$, the vector field $V_{\alpha}:=V_{\alpha_+}:=\left({proj}_{\Lambda,+}^g \alpha \right)^{\sharp}$, which is independent of $g \in c$, is a normal conformal vector field, i.e. the dual of a nc-Killing 1-form. We denote the space of all normal conformal vector fields on $(M,g)$ by $\mathfrak{X}^{nc}(M)$. \cite{raj} shows that for a vector field $V$ being normal conformal is equivalent to being conformal, $V \in \mathfrak{X}^c(M)$, i.e. $L_Vg = \lambda \cdot g$, and to satisfy in addition that
\begin{equation}
\begin{aligned}
V \invneg W^g = 0,\text{ }V \invneg C^g = 0,
\end{aligned}
\end{equation}
where $W^g$ and $C^g$ are the Weyl- and Cotton-York tensor, respectively.

\section{The general construction of tractor conformal superalgebras} \label{cts}
Let $(M^{1,n-1},c)$ be a connected, oriented and time-oriented Lorentzian conformal spin manifold. Here, when dealing with spinor- and spin tractor bundles, we always mean the \textit{complex ones}, i.e. $\mathcal{S}(M)$ or $S^g(M)$ are obtained as associated vector bundles to $\mathcal{Q}^1$ or $\mathcal{Q}^g$ using $\Delta_{2,n}^{\C}$ or $\Delta^{\C}_{1,n-1}$, respectively.

The previous section revealed that distinguished classical and spinorial conformal symmetries of $(M,c)$ can be described in a conformally invariant way in terms of parallel tractors. For the construction of a superalgebra canonically associated to $M,c)$, we therefore set
\begin{align*}
\g_0 &:= {Par}\left(\Lambda_{\mathcal{T}}^2(M), \nabla^{nc}\right) \subset \Omega^2_{\mathcal{T}}(M).\\
\g_1 &:= {Par}\left(\mathcal{S}(M), \nabla^{nc}\right) \subset \Gamma\left(\mathcal{S}(M), \nabla^{nc}\right).
\end{align*}
By means of $g \in c$, $\g_0$ and $\g_1$ correspond to normal conformal vector fields and twistor spinors, respectively. We now introduce natural brackets which make $\g=\g_0 \oplus \g_1$ become a superalgebra:\\
\newline
For the even-even bracket, we first globalize the isomorphism $\mathfrak{so}(2,n) \cong \Lambda^2_{2,n}$ to obtain
\begin{equation}\label{suz}
\begin{aligned}
\tau : \Omega^2_{\mathcal{T}}(M) &\rightarrow  \mathfrak{so}(\mathcal{T}(M),\langle \cdot, \cdot \rangle_{\mathcal{T}}),\\
\alpha & \mapsto  \alpha_E \text {, }\alpha_E(X):=(X \invneg \alpha)^{\sharp}. 
\end{aligned}
\end{equation}
$ \mathfrak{so}(\mathcal{T}(M),\langle \cdot, \cdot \rangle_{\mathcal{T}})$ carries the pointwise defined usual Lie bracket of endomorphisms. We set for $\alpha, \beta \in \g_0$
\begin{align*}
[\alpha, \beta]:= \tau^{-1}\left( \alpha_E \circ \beta_E - \beta_E \circ \alpha_E \right).
\end{align*}
Moreover, $\nabla^{nc}$ induces a covariant derivative $\nabla^{nc}$ on $ \mathfrak{so}(\mathcal{T}(M),\langle \cdot, \cdot \rangle_{\mathcal{T}})$ in a natural way.
\begin{Proposition}
For $\alpha,\beta \in \g_0$ we have that also $[\alpha,\beta] \in \g_0$.
\end{Proposition}
\textit{Proof. }We first show that that $\alpha \in {Par}\left(\Lambda_{\mathcal{T}}^2(M),\nabla^{nc} \right) \Leftrightarrow \alpha_E \in {Par}\left(\mathfrak{so}(\mathcal{T}(M),\langle \cdot, \cdot \rangle_{\mathcal{T}}),\nabla^{nc} \right)$: Let $X \in \mathfrak{X}(M)$, $x \in M$ and let $(v_0,...,v_{n+1})$ be a local frame in $\mathcal{T}M$ which is parallel in $x$ wrt. $\nabla^{nc}$. We have for $i \in \{0,...,n+1 \}$ at $x$:
\begin{align*}
\left(\nabla^{nc}_X \alpha_E \right)(v_i) & = \nabla^{nc}_X \left( \alpha_E (v_i) \right)= \nabla^{nc}_X {\left( v_i \invneg \alpha \right)}^{\sharp} = \left( \nabla^{nc}_X \left( v_i \invneg \alpha \right)\right)^{\sharp} = (v_i \invneg \nabla^{nc}_X \alpha)^{\sharp},
\end{align*}
which proves this claim. Thus, it suffices to check that for $\alpha, \beta \in \g_0$ also $[\alpha_E, \beta_E]_{\mathfrak{so}} \in {Par}\left(\mathfrak{so}(\mathcal{T}(M),\langle \cdot, \cdot \rangle_{\mathcal{T}}),\nabla^{nc} \right)$. We compute with the same notations as above at $x$:
\begin{align*}
\left(\nabla_X^{nc} \left([\alpha_E, \beta_E]_{\mathfrak{so}} \right) \right)(v_i) &= \nabla_X^{nc} \left([\alpha_E, \beta_E]_{\mathfrak{so}}(v_i) \right) - [\alpha_E, \beta_E]_{\mathfrak{so}} \left( \nabla_X^{nc} v_i \right) \\
&=\nabla_X^{nc} \left( \alpha_E (\beta_E (v_i))-\beta_E (\alpha_E(v_i)) \right) =  \alpha_E (\beta_E (\nabla_X^{nc}v_i))-\beta_E (\alpha_E (\nabla_X^{nc}v_i))\\
&=0
\end{align*}
This proves the Proposition.
$\hfill \Box$\\

Clearly, $\g_0$ now becomes a Lie algebra in the usual sense. We shall show in the next section that the chosen bracket is \textit{the right one} in the sense that if $\alpha, \beta$ are considered as normal conformal vector fields for some fixed $g \in c$ by means of $\left( {proj}^g_{\Lambda,+} \alpha \right)^{\sharp}$, then $[ \cdot, \cdot ]$ translates into the usual Lie bracket of vector fields. \\
\newline
As a next step we define the odd-odd bracket, which by definition has to be a symmetric bilinear map $\g_1 \times \g_1 \rightarrow \g_0$. A nontrivial way to obtain a parallel tractor 2-form from two parallel spin tractors is given by the parallel tractor form (\ref{3245}), i.e.
\begin{align*}
[ \cdot , \cdot ] : \g_1 \times \g_1 \rightarrow \g_0 \text{ , } (\psi_1,\psi_2) \mapsto \alpha_{\psi_1,\psi_2}^2.
\end{align*}
In signature $(2,n)$, the form $\alpha_{\psi_1,\psi_2}^2$ is given as follows: One observes that $\langle \alpha \cdot \psi, \psi \rangle_{\Delta_{2,n}^{\C}} \in i\R$ for $\psi \in \Delta_{2,n}^{\C}, \alpha \in \Lambda^2_{2,n}$. (\ref{3245}) thus yields that
\begin{align}
\langle \alpha_{\psi_1,\psi_2}^2, \alpha \rangle_{\mathcal{T}} = \text{Im }\langle \alpha \cdot \psi_1, \psi_2 \rangle_{\mathcal{S}}\text{,  }\alpha \in \Omega^2_{\mathcal{T}}(M). \label{impa}
\end{align}
$\alpha_{\psi_1,\psi_2}^2$ is then symmetric in $\psi_1$ and $\psi_2$.\\
\newline
It remains to introduce an even-odd-bracket. We set
\begin{align*}
[ \cdot , \cdot ] : \g_0 \times \g_1 \rightarrow \g_1 \text{ , } (\alpha ,\psi) \mapsto \frac{1}{2} \alpha \cdot \psi.
\end{align*}
The meaning of the factor $\frac{1}{2}$ will become clear in a moment. It follows directly from $\nabla_X^{nc} (\alpha \cdot \psi) = (\nabla_X^{nc}\alpha) \cdot \psi + \alpha \cdot \nabla_X^{nc} \psi$ that this map is well-defined, i.e. the image lies again in $\g_1$.  Moreover, in order to obtain the right symmetry relations, we must set 
\begin{align*}
[ \cdot , \cdot ] : \g_1 \times \g_0 \rightarrow \g_1 \text{ , } (\psi ,\alpha) \mapsto -\frac{1}{2} \alpha \cdot \psi.
\end{align*}
With these choices of $\g_0, \g_1$ and definitions of the brackets, we have associated a nontrivial (real) conformal superalgebra to the conformal structure (where $\g_1$ is considered as a \textit{real} vector space). 

\begin{definition}
The (real) superalgebra $\g=\g_0 \oplus \g_1$ associated to $(M^{1,n-1},c)$ is called the tractor conformal superalgebra (associated to $(M,c)$).
\end{definition}

It is natural to ask under which circumstances the construction produces a \textit{Lie} superalgebra, i.e. we have to check the four Jacobi identities from (\ref{grj}). As $\g_0$ is a Lie algebra in its own right, the even-even-even Jacobi identity is always satisfied.

\begin{Proposition}
The tractor conformal superalgebra associated to a Lorentzian conformal spin manifold satisfies the even-even-odd and the even-odd-odd Jacobi identity.
\end{Proposition}

\textit{Proof. }By (\ref{grj}) we have to check that 
\begin{align*}
[\alpha,[\beta,\psi]]\stackrel{!}{=}[[\alpha,\beta],\psi]+[\beta,[\alpha,\psi]] \text{ }\forall \alpha,\beta \in \g_0, \psi \in \g_1,
\end{align*}
which by definition of the brackets is equivalent to showing that
\begin{align*}
2 \cdot [\alpha,\beta] \cdot \psi \stackrel{!}{=} \alpha \cdot \beta \cdot \psi - \beta \cdot \alpha \cdot \psi,
\end{align*}
being a purely algebraic identity at each point. Whence, we may for the proof assume that $\alpha, \beta \in \Lambda^2_{2,n}$ and $\psi \in \Delta_{2,n}^{\C}$. With respect to the standard basis of $\R^{2,n}$ we express
\begin{align*}
\alpha = \sum_{i<j} \epsilon_i \epsilon_j \alpha_{ij} e^{\flat}_i \wedge e^{\flat}_j \Rightarrow \alpha_E = \sum_{i<j} \epsilon_i \epsilon_j \alpha_{ij} E_{ij} \text{ and }
\beta = \sum_{k<l} \epsilon_k \epsilon_l \beta_{kl} e^{\flat}_k \wedge e^{\flat}_l \Rightarrow \beta_E = \sum_{k<l} \epsilon_k \epsilon_l \beta_{kl} E_{kl}.
\end{align*}
Here, $E_{kl}:= \epsilon_k D_{lk} - \epsilon_l D_{kl}$  with $k<l$ form a basis of the Lie algebra $\mathfrak{so}(2,n)$, where $D_{kl}$ denotes the matrix in $M(n+2,\R)=\mathfrak{gl}(n+2,\R)$ whose $(k,l)$ entry is 1 and all other entries are 0. The Lie algebra relations read
\begin{align}
[E_{ij},E_{kl}]_{\mathfrak{so}(2,n)}=\begin{cases} 0 & i=k, j=l\text{ or }i,j,k,l \text{ pairwise distinct,} \\ \epsilon_i E_{jl} & i=k, j \neq l, \\ 
\end{cases} \label{so}
\end{align}
This shows that 
{\allowdisplaybreaks
\begin{align*}
2 \cdot [\alpha, \beta] \cdot \psi & = \tau^{-1}\left([\alpha_E,\beta_E]_{\mathfrak{so}(2,n)}\right) = \sum_{i<j}\sum_{k<l} \epsilon_i \epsilon_j \epsilon_k \epsilon_l \alpha_{ij} \beta_{kl} \tau^{-1}\left(2 \cdot [E_{ij},E_{kl}]_{\mathfrak{so}(2,n)}\right) \cdot \psi \\ 
& =  \sum_{i<j}\sum_{k<l} \epsilon_i \epsilon_j \epsilon_k \epsilon_l \alpha_{ij} \beta_{kl} [e_ie_j,e_ke_l]_{\mathfrak{spin}(2,n)} \cdot \psi = (\alpha \cdot \beta - \beta \cdot \alpha) \cdot \psi.
\end{align*}}
The even-odd-odd Jacobi identity is by polarization equivalent to $[\alpha,[\psi,\psi]]=[[\alpha,\psi],\psi]+[\psi,[\alpha,\psi]]$ for all $\alpha \in \g_0$ and $\psi \in \g_1$. By definition of the brackets, we have to show that
\begin{align}
\left[\alpha_E,\left(\alpha_{\psi}^2\right)_E\right]_{\mathfrak{so}(\mathcal{T}(M))}\stackrel{!}{=}\left(\frac{1}{2}\alpha^2_{\alpha \cdot \psi, \psi} + \frac{1}{2}\alpha^2_{\psi,\alpha \cdot  \psi}\right)_E = \left(\alpha^2_{\alpha \cdot \psi, \psi}\right)_E. \label{67}
\end{align}
Again, this is pointwise a purely algebraic identity. Whence, it suffices to prove it for $\alpha \in \Lambda^2_{2,n}$ and $\psi \in \Delta^{\C}_{2,n}$. With respect to the standard basis of $\R^{2,n}$, we write $\alpha$ and $\alpha_E$ as above. Inserting the definition of $\alpha_{\psi}^2$ leads to
\begin{align} \label{w1}
\left[\alpha_E,\left(\alpha_{\psi}^2\right)_E \right]= \sum_{i<j}\sum_{k<l} \epsilon_i \epsilon_j \epsilon_k \epsilon_l \alpha_{ij} \cdot \text{Im }(\langle e_k \cdot e_l \cdot \psi, \psi \rangle_{\Delta_{2,n}^{\C}}) \cdot [E_{ij},E_{kl}],
\end{align}
whereas the right-hand side of (\ref{67}) is by definition given by
\begin{align} 
\left(\alpha^2_{\alpha \cdot \psi, \psi}\right)_E & = \sum_{k < l} \epsilon_k \epsilon_l \text{Im }(\langle e_k \cdot e_l \cdot \alpha \cdot \psi, \psi \rangle_{\Delta_{2,n}^{\C}}) \cdot E_{kl} \notag \\
& = \sum_{i < j} \sum_{k < l} \epsilon_i \epsilon_j \epsilon_k \epsilon_l \alpha_{ij} \cdot \text{Im }(\langle e_k \cdot e_l \cdot e_i \cdot e_j \cdot \psi, \psi \rangle_{\Delta_{2,n}^{\C}}) \cdot E_{kl}. \label{w2}
\end{align}
Using the algebra relations for $\mathfrak{so}(2,n)$ from (\ref{so}), it is not difficult to show that every summand in (\ref{w1}) shows up also in (\ref{w2}) and vice versa:
\begin{itemize}
\item Consider summands with $i,j,k,l$ pairwise distinct. Clearly, they vanish in (\ref{w1}). On the other hand, $\langle e_k \cdot e_l \cdot e_i \cdot e_j \cdot \psi, \psi \rangle_{\Delta_{2,n}^{\C}} \in \R$, i.e. the summands also vanish in (\ref{w2}).
\item Consider summands with $i=k,j=l$. Again, they vanish in (\ref{w1}). In (\ref{w2}), these summands are proportional to $\langle \psi,\psi \rangle_{\Delta^{\C}_{2,n}} \in \R$, so the imaginary part vanishes.
\item Consider summands in (\ref{w2}) with $i=k$ and $j \neq l$. They lead to the expression $ -\epsilon_j \epsilon_l \alpha_{ij} \text{Im }(\langle \epsilon_i \cdot e_j \cdot e_l \cdot \psi, \psi \rangle_{\Delta^{\C}_{2,n}}) E_{il}$. In (\ref{w1}), these summands can be found for choosing $j=k$ and $i\neq l$ for which we get $[E_{ij},E_{kl}]=-\epsilon_j E_{il}$, and thus the summand $ -\epsilon_j \epsilon_l \alpha_{ij} \text{Im }(\langle \epsilon_i \cdot e_j \cdot e_l \cdot \psi, \psi \rangle_{\Delta^{\C}_{2,n}}) E_{il}$ also shows up in (\ref{w1}). The remaining cases are equivalent to this one after permuting the indices.
\end{itemize}
Consequently, the two sums are identical and (\ref{67}) holds. 
$\hfill \Box$\\

In contrast to that, the remaining Jacobi identity does not hold in general as we shall later see for concrete examples. Under certain restrictions on the conformal holonomy representation, we can however show that all Jacobi identities hold.

\begin{satz} \label{hola}
Suppose that the conformal holonomy representation of $Hol_x(M,c)$ on $\mathcal{T}_x(M)$ for $x \in M$ satisfies the following: There exists \textbf{no} (possibly trivial) $m-$dimensional Euclidean subspace $E \subset \mathcal{T}_x(M) \cong \R^{2,n}$ such that both
\begin{enumerate}
\item The action of $Hol_x(M,c)$ fixes $E$ (and therefore also $E^{\bot}$).
\item $E^{\bot}$ is even-dimensional and on $E^{\bot}$, $Hol_x(M,c)_{E^{\bot}}:=\{A_{|E^{\bot}} \mid A \in Hol_x(M,c) \} \subset SO^+(E^{\bot})\cong SO^+(2,n-m)$ is conjugate to a subgroup of $SU(1,\frac{n-m}{2})\subset SO(2,n-m)$.
\end{enumerate}
Then the tractor conformal superalgebra $\g$ satisfies the odd-odd-odd Jacobi identity, and thus carries the structure of a Lie superalgebra.
\end{satz}

\textit{Proof. }As a first step, we show that under the assumptions, 
\[\psi \in \g_1 \Rightarrow \text{ker }\psi:= \{ v \in \mathcal{T}(M) \mid v \cdot \psi = 0 \} \neq \{0 \}.\] 

To this end, note that all possible algebraic Dirac forms $\alpha_{\chi}^2$ for $0 \neq \chi \in \Delta^{\C}_{2,n-2}$ have been classified in \cite{leihabil}. Precisely one of the following cases occurs:
\begin{enumerate}
\item $\alpha_{\chi}^2 = l_1^{\flat} \wedge l_2^{\flat}$, where $l_1,l_2$ span a totally lightlike plane in $\R^{2,n-2}$.
\item $\alpha_{\chi}^2 = l^{\flat} \wedge t^{\flat}$ where $l$ is lightlike, $t$ is a orthogonal timelike vector.
\item $\alpha_{\chi}^2 =\omega$ (up to conjugation in $SO(2,n-2)$), where $n=2m$ is even and $\omega$ is equivalent to the standard K\"ahler form $\omega_0$ \footnote{By this we mean that there are nonzero constants $\mu_i \in \R$ such that $\omega = \sum_{i=1}^m \mu_i e_{2i-1}^{\flat} \wedge  e_{2i}^{\flat}$. One obtains $\omega_0$, the standard pseudo-K\"ahler form for $\mu_i=1$ for all $i$.} on $\R^{2,n-2}$. In this case, $Stab_{\alpha_{\chi}^2} O(2,n-2) \subset U(1,m-1)$.
\item There is a nontrivial Euclidean subspace $E \subset \R^{2,n-2}$ such that ${\alpha_{\chi}^2}_{|E} = 0$ and $\alpha_{\chi}^2$ is equivalent to the standard K\"ahler form on the orthogonal complement $E^{\bot}$ of signature $(2,2m)$ (again, this is up to conjugation in $SO(2,n-2)$). In this case $Stab_{\alpha_{\chi}^2} O(2,n-2) \subset U(1,m) \times O(n-2(m+1))$.
\end{enumerate}
Moreover, one easily calculates that the first case occurs iff ker $ \chi$ is 2-dimensional (and in this case it is spanned by $l_1,l_2$). The second case occurs iff this kernel is one-dimensional (and spanned by $l$), whereas the last two cases can only occur if the kernel under Clifford multiplication is trivial. \\
For $\psi \in \g_1$, the parallel tractor 2-form $\alpha_{\psi}^2$, whose $SO^+(2,n)$-orbit type is constant over $M$, must up to conjugation be one of the four generic types from the above list. Types 3. and 4. obviously contradict our assumptions. Whence, $\alpha_{\psi}^2$ is of type 1. or 2, yielding that $\text{dim ker }\psi \in \{1,2\}$.\\
By a standard polarization argument the odd-odd-odd Jacobi identity is equivalent to show that $[\psi,[\psi,\psi]]=0$ for all $\psi \in \g_1$. By definition of the brackets, this precisely says that
\begin{align*}
\alpha_{\psi}^2 \cdot \psi \stackrel{!}{=} 0.
\end{align*}
However, as ker $\psi \neq \{0 \}$, the above discussion yields that $\alpha_{\psi}^2 = l^{\flat} \wedge r^{\flat}$, where $l \in \text{ker }\psi$ and $r$ is orthogonal to $l$. It follows that $\alpha_{\psi}^2 \cdot \psi = - r \cdot l \cdot \psi = 0$. This proves the remaining Jacobi identity and the Theorem.
$\hfill \Box$ \\
\newline 
The requirements from Theorem \ref{hola} translate into more down-to-earth geometric statements using the classification of Lorentzian manifolds admitting twistor spinors by F. Leitner:

\begin{satz}[\cite{leihabil}; Thm.10] \label{bg} 
Let $\varphi=\widetilde{\Phi}^g ( {proj}^g_+ \psi) \in \Gamma(S^g)$ be a complex twistor spinor on a Lorentzian spin manifold $(M^{1,n-1},g)$ of dimension $n\geq 3$. Then one of the following holds on an open and dense subset $\widetilde{M} \subset M$:
\begin{enumerate}
\item $\alpha_{\psi}^2=l_1^{\flat} \wedge l_2^{\flat}$ for standard tractors $l_1,l_2$ which span a totally lightlike plane.\\
In this case, $\ph$ is locally conformally equivalent to a parallel spinor with lightlike Dirac current $V_{\ph}$ on a Brinkmann space.
\item $\alpha_{\psi}^2=l^{\flat} \wedge t^{\flat}$ where $l$ is a lightlike, $t$ is an orthogonal, timelike standard tractor.\\
$(M,g)$ is locally conformally equivalent to $(\mathbb{R},-dt^2) \times (N_1, h_1) \times \cdots \times (N_r, h_r)$, where the $(N_i,h_i)$ are Ricci-flat K\"ahler, hyper-K\"ahler, $G_2$-or $Spin(7)$-manifolds.
\item $\alpha_{\psi}^2$ is of K\"ahler-type (cf. cases 3. and 4. in the above list)\\
The following cases can occur:
\begin{enumerate}
\item The dimension $n$ is odd and the space is locally equivalent to a Lorentzian Einstein-Sasaki manifold on which the spinor is a sum of Killing spinors.
\item $n$ is even and $(M,g)$ is locally conformally equivalent to a Fefferman space.
\item 
There exists locally a product metric $g_1 \times g_2 \in [g]$ on $M$, where $g_1$ is a Lorentzian Einstein-Sasaki metric on a space $M_1$ of dimension $n_1=2 \cdot rk(\alpha_1(\ph))+1 \geq 3$ admitting a Killing spinor and $g_2$ is a Riemannian Einstein metric with Killing spinor on a space $M_2$ of positive scalar curvature $scal^{g_2} = \frac{(n-n_1)(n-n_1-1)}{n_1 (n_1 - 1)}scal^{g_1}$.
\end{enumerate}
\end{enumerate} 
\end{satz}
Again, ker ${\psi} \neq \{0 \}$ occurs exactly in the first two cases of Theorem \ref{bg}. In the third case of Theorem \ref{bg}, it hold that dim ker ${\psi}=\{0\}$
Consequently, geometries admitting twistor spinors and which do \textbf{not} satisfy the conditions from Theorem \ref{hola} correspond to the cases $3.(a)-3.(c)$ mentioned in Theorem \ref{bg}, being Fefferman metrics, Lorentzian Einstein Sasaki manifolds or local splittings $g_1 \times g_2 \in [g]$ where $g_1$ is a Lorentzian Einstein Sasaki metric and $g_2$ is a Riemannian Einstein metric of positive scalar curvature.  Thus, Theorem \ref{hola} can be rephrased in more geometric terms by saying that if none of these three special geometries lies in the conformal class of the metric, one obtains a conformal tractor Lie superalgebra. \\
This is in accordance with other observations in the literature (cf. \cite{med}). Namely it is known that for the mentioned special geometries one has to include further symmetries in the algebra in order to obtain a conformal Lie superalgebra to which we will come back later.

\begin{bemerkung}
The construction of a \textit{real} tractor conformal superalgebra can completely analogous be carried out with \textit{real} spinors. One then has to make the obvious modifications, i.e. define $\alpha^2_{\psi_1,\psi_2}$ without the imaginary part from (\ref{impa}). Note that $\langle \psi, \psi \rangle_{\Delta_{2,n}^{\R}} = 0 \text{ }\forall \psi \in \Delta_{2,n}^{\R}$. One obtains the same results, i.e. all Jacobi identities except the odd-odd-odd one are always satisfied. However, as we are later dealing with tractor conformal superalgebras for twistor spinors on Fefferman spaces, cf. \cite{bafe}, it seems more appropriate to work with complex quantities in this chapter.
\end{bemerkung}

\begin{bemerkung} \label{remak}
We defined the even part of the tractor conformal superalgebra to be (isomorphic to) the space of \textit{normal} conformal vector fields. It is possible to include all conformal vector fields $\mathfrak{X}^c(M)$ in the even part using tractor calculus as follows: Let $\alpha \in \Omega^2_{\mathcal{T}}(M)$ be a tractor 2-form on $(M,c)$ and let $V_{\alpha}=\left( {proj}^g_{\Lambda,+}(\alpha) \right)^{\sharp} \in \mathfrak{X}(M)$ be the associated vector field. As proved in \cite{cov1,cov2}, we have that $V_{\alpha} \in \mathfrak{X}^c(M)$, the space of conformal vector fields, if and only if
\begin{align}
\nabla_X^{nc} \alpha = \tau^{-1} \left({R}^{\nabla^{nc},\mathcal{T}(M)}(V_{\alpha},X) \right) \text{ }\forall X \in \mathfrak{X}(M), \label{star}
\end{align}
where we identify the skew-symmetric curvature endomorphism with a tractor 2-form by means of the isomorphism $\tau$ from (\ref{suz}). We now consider the extended tractor superalgebra
\begin{align*}
\mathfrak{g}_0^{ec}:= \{ \alpha \in \Omega^2_{\mathcal{T}}(M) \mid \alpha \text{ satisfies } (\ref{star}) \} \text{ and } \mathfrak{g}^{ec}:=\mathfrak{g}_0^{ec} \oplus \g_1,
\end{align*}
where $\g_1 = Par(\mathcal{S},\nabla^{nc})$ is as before. On this space, we may define the same brackets as defined on $\mathfrak{g}$ above and observe that they are still well-defined: For $\alpha,\beta \in \mathfrak{g}_0^{ec}$, we have that also $[\alpha,\beta] \in \mathfrak{g}_0^{ec}$ as by Proposition \ref{cof} $V_{[\alpha,\beta]} = -[V_{\alpha},V_{\beta}]_{\mathfrak{X}(M)}$, which is a conformal vector field. Next, let $\alpha \in \mathfrak{g}_0^{ec}$ and $\psi \in \g_1$. Then we have that
\begin{align*}
\nabla^{nc}_X (\alpha \cdot \psi) &=  (\nabla_X^{nc} \alpha ) \cdot \psi = \tau^{-1} \left({R}^{\nabla^{nc},\mathcal{T}(M)}(V_{\alpha},X) \right) \cdot \psi \\
&= 2 \cdot {R}^{\nabla^{nc},\mathcal{S}}(V_{\alpha},X)\psi \stackrel{\psi \in \mathfrak{g}_1}{=} 0,
\end{align*}
i.e. $\alpha \cdot \psi \in \mathfrak{g}_1.$ This shows that $\mathfrak{g}^{ec}$ together with the defined brackets is a conformal superalgebra which naturally extends $\mathfrak{g}$. Moreover, the subsequent Propositions \ref{cof} and \ref{sld} and (\ref{765}) still hold in this situation and describe $\mathfrak{g}^{ec}$ wrt. a metric $g \in c$ as their proofs only involve the conformal Killing equation for vector fields and not the normalisation conditions. \\
However, we will only consider the superalgebra $\mathfrak{g}$ and not its extension $\mathfrak{g}^{ec}$ in the sequel because in the case of twistor spinors there are always normal conformal vector fields, and it seems to us that the structure of the subalgebra $\mathfrak{g}$ and the existence of distinguished normal conformal vector fields is more directly related to special geometric structures (cf. \cite{nc}) on $(M,c)$ than the structure of $\mathfrak{g}^{ec}$ as we will see in the next sections.

\end{bemerkung}

\section{Metric description and examples} \label{dem}
Fixing a metric $g \in c$ leads to canonical isomorphisms
\begin{align*}
i_0 &: \g_0 = {Par}\left(\Lambda_{\mathcal{T}}^2(M), \nabla^{nc}\right) \rightarrow  \mathfrak{X}^{nc}(M), &\alpha &\mapsto V_{\alpha}:=\left({proj}^g_{\Lambda,+}(\alpha)\right)^{\sharp}, \\
i_1 &: \g_1 = {Par}\left(\mathcal{S}(M), \nabla^{nc}\right) \rightarrow \text{ker }P^g, &\psi &\mapsto \ph:=\widetilde{\Phi}^g({proj}^g_+(\psi)).
\end{align*}
The aim of this section is to compute the behaviour of the tractor conformal superalgebra structure under these isomorphisms. As it turns out, the maps $i_0$ and $i_1$ allow us to identify our tractor conformal superalgebra with conformal superalgebras constructed for Lorentzian conformal spin manifolds in \cite{raj} and \cite{ha96}.

\begin{Proposition} \label{cof}
For fixed $g \in c$ it holds for all $\alpha, \beta \in \g_0$ that
\begin{align*}
i_0\left([\alpha, \beta ]_{\g_0}\right) = \left[V_{\beta},V_{\alpha}\right]_{\mathfrak{X}(M)} = \left[i_0(\beta),i_0(\alpha)\right]_{\mathfrak{X}(M)}
\end{align*}
\end{Proposition}
\textit{Proof. }
We start with some algebraic computations: Assume that $\alpha, \beta \in \Lambda^2_{2,n}$. Wrt. the decomposition (\ref{deco}) we may write
$\alpha = e_+^{\flat} \wedge \alpha_+ + \alpha_0 + \alpha_{\mp} \cdot e_-^{\flat} \wedge e_+^{\flat} + e_-^{\flat} \wedge \alpha_-$ with $\alpha_+ = \sum_{i=1}^{n}\epsilon_i \alpha_i^+ \cdot e_i^{\flat}$, $\alpha_-=\sum_{i=1}^{n}\epsilon_i \alpha_i^- \cdot e_i^{\flat}$, $\alpha_{0} = \sum_{i<j} \epsilon_i \epsilon_j \alpha_{ij}^0 \cdot e_i^{\flat} \wedge e_j^{\flat}$ for real coefficients $\alpha_i^+$ etc. For the standard basis of the Lie algebra $\mathfrak{so}(2,n)$, cf. \ref{so}, we let $E_{\pm, i} := \frac{1}{\sqrt{2}}(E_{n+1 i}\pm E_{0i})$. Then the endomorphism $\alpha_E = \tau(\alpha) \in \mathfrak{so}(2,n)$ is given by
\begin{align*}
\alpha_E = \sum^n_{i=1} \epsilon_i \alpha_i^+ E_{+i} + \sum^n_{i=1} \epsilon_i \alpha_i^- E_{-i} + \alpha_{\mp} E_{n+1 0} + \sum_{i<j}^n \epsilon_i \epsilon_j \alpha_{ij}^0 E_{ij}.
\end{align*}
An analogous expression holds for $\beta_E$. Using the algebra relations (\ref{so}), it is straightforward to compute the following commutators for $i,j=1,...,n$:
\begin{align*}
[E_{\pm,i},E_{\pm,j}]&=0, \\
[E_{-,i},E_{+,j}]&=E_{ij}-\epsilon_i \delta_{ij} E_{0n+1}+\epsilon_j \delta_{ij} E_{0n+1}, \\
[E_{\pm,i},E_{n+1 0}]&=\mp E_{\pm,i}, \\
[E_{ij},E_{\pm,k}]&=\epsilon_i \delta_{ik} E_{\pm,j} - \epsilon_j \delta_{jk} E_{\pm,i}.
\end{align*}
With these formulas, we compute
\begin{align*}
\left[\alpha_E,\beta_E \right]_{\mathfrak{so}(2,n)} =&+ \sum_{i=1}^n \epsilon_i (\beta_i^+ \alpha_{\mp} - \alpha_i^+ \beta_{\mp} ) E_{+,i} + \sum_{i<j} \epsilon_i \epsilon_j (\alpha_{ij}^0 \beta_i^+ - \beta_{ij}^0 \alpha_i^+) E_{+,j} \\
&- \sum_{j<i} \epsilon_i \epsilon_j (\alpha_{ji}^0 \beta_j^+ - \beta_{ji}^0 \alpha_j^+) E_{+,i} + \text{Terms not involving }E_{+,i}.
\end{align*}
A global version of this formula yields that for $\alpha,\beta \in \g_0$ one has wrt. $g \in c$
\begin{align*}
{proj}_{\Lambda,+}^g\left([\alpha,\beta]_{\g_0}\right) = \alpha_{\mp} \cdot \beta_+ - \beta_{\mp} \alpha_+ + {\sum_{i < j} \epsilon_i \epsilon_j (\alpha_{ij}^0 \beta_i^+  - \beta_{ij}^0  \alpha_i^+) s_j^{\flat}}- \sum_{j<i} \epsilon_i \epsilon_j (\alpha_{ji}^0 \beta_j^+ - \beta_{ji}^0 \alpha_j^+) s_i^{\flat},
\end{align*}
where $(s_1,...,s_n)$ is a local $g-$pseudo-orthonormal frame in $TM$, with coefficients of $\alpha$ taken with respect to this frame. This can be rewritten as 
\begin{align}
i_0\left([\alpha,\beta]_{\g_0}\right) =\left({proj}_{\Lambda,+}^g\left([\alpha,\beta]_{\g_0}\right)\right)^{\sharp} = \alpha_{\mp} \cdot V_{\beta} - \beta_{\mp} V_{\alpha} + \left(V_{\beta} \invneg \alpha_0 - V_{\alpha} \invneg \beta_0 \right)^{\sharp} \label{e1}
\end{align}
We now compare this expression to the Lie bracket $[V_{\alpha},V_{\beta}]$. Dualizing the first nc-Killing equation (cf. \ref{stu}) for $\alpha_+$ yields that
\begin{align*}
\nabla^g_X V_{\alpha} = (X \invneg \alpha_0)^{\sharp} + \alpha_{\mp} \cdot X \text{ }\forall X \in \mathfrak{X}(M).
\end{align*}
Consequently, 
\begin{align}
[V_{\beta},V_{\alpha}]=\nabla^g_{V_{\beta}} V_{\alpha}-\nabla^g_{V_{\alpha}}V_{\beta} = \left(\alpha_{\mp}V_{\beta} - \beta_{\mp} V_{\alpha}\right) + \left(V_{\beta} \invneg \alpha_0 - V_{\alpha} \invneg \beta_0 \right)^{\sharp}. \label{e2}
\end{align}
Comparing the two expressions (\ref{e1}) and (\ref{e2}) immediately yields the claim.
$\hfill \Box$\\

The next Proposition will be proved in a more general setting in Proposition \ref{lsd}:

\begin{Proposition} \label{sld}
For $\alpha \in \g_0$, $\psi \in \g_1$, and $g \in c$ such that $\ph = \widetilde{\Phi}^g\left({proj}^g_+ \psi \right)=i_1(\psi)$, and $V_{\alpha_+} = i_0(\alpha)$, we have that 
\begin{align*}
i_1\left(\left[\alpha,\psi \right]_{\g_1}\right) = \frac{1}{2}(\widetilde{\Phi}^g \circ {proj}_+^g) \left(\alpha \cdot \psi \right) = -  \underbrace{\left(\nabla_{V_{\alpha}} \ph + \frac{1}{4} \tau \left(\nabla V_{\alpha} \right) \cdot \ph \right) }_{=:V_{\alpha} \circ \ph}, 
\end{align*}
where $\tau \left(\nabla V_{\alpha} \right) := \sum_{j=1}^n \epsilon_j \left( \nabla_{s_j} V_{\alpha} \right) \cdot s_j + (n-2) \cdot \lambda_{\alpha_+}$ and $L_{V_{\alpha_+}}g = 2 \lambda_{\alpha_+}g$. 
\end{Proposition}

\begin{bemerkung}
The above term $V_{\alpha} \circ \ph$ is the spinorial Lie derivative used in \cite{kos,ha96,raj} for the construction of a conformal Killing superalgebra.
\end{bemerkung}

Finally, we give the metric expression of the odd-odd bracket. Let $\psi \in \g_1$ and $\ph=\widetilde{\Phi}^g \left({proj}^g_+(\psi)\right)=i_1(\psi)$:
\begin{align}
i_1 \left( [\psi,\psi] \right) &= \left( {proj}^g_{\Lambda,+} \left( \alpha_{\psi}^2 \right) \right)^{\sharp} \stackrel{(\ref{tuffi})}{=} c_{1,1}^1 \cdot \left( \alpha^1_{\ph} \right)^{\sharp} = c_{1,1}^1 \cdot V_{\ph}, \label{765}
\end{align}
where the nonzero constant $c_{1,1}^1 \in \R$ from (\ref{tuffi}) depends only on the choice of an admissible scalar product (in the sense of \cite{cor}) on $\Delta_{2,n}^{\C}$. These computations directly prove the following statement:

\begin{satz}
Given a Lorentzian conformal spin manifold $(M^{1,n-1},c)$, the associated tractor conformal superalgebra $\g=\g_0 \oplus \g_1 = {Par}\left(\Lambda_{\mathcal{T}}^2(M), \nabla^{nc}\right) \oplus {Par}\left(\mathcal{S}(M), \nabla^{nc}\right)$ is via a fixed $g \in c$ isomorphic to the conformal superalgebra (\ref{alr}) on $\mathfrak{X}^{nc}(M)\oplus \text{ker }P^g$ (as considered in \cite{raj}). Up to prefactors, the $g-$dependent maps $i_0$ and $i_1$ are superalgebra (anti-)isomorphisms.
\end{satz}

\begin{bemerkung}
If for some fixed $g \in c$ the manifold admits geometric Killing spinors, for instance if there is an Einstein metric in the conformal class, the restrictions of the brackets (\ref{alr}) to $\mathfrak{X}^k(M) \oplus \mathcal{K}(M)$, the space of Killing vector fields and Killing spinors as even and odd parts, is well-defined (cf. \cite{suads}), and thus gives a subalgebra of the superalgebra $\mathfrak{g}^{ec}$.\\ 
More generally, the construction of Killing superalgebras for Riemannian or Lorentzian manifolds using the cone construction where the even part consists of Killing vector fields and the odd part of geometric Killing spinors is discussed in \cite{suads}. In case of an Einstein metric in the conformal class this is equivalent to our tractor construction as in this case all conformal holonomy computations restrict to considerations on the metric cone, see \cite{leihabil,baju}.
\end{bemerkung}

Let us consider some examples of tractor conformal superalgebras:
\subsection*{Tractor conformal superalgebras with one twistor spinor}
Consider the case that dim $\g_1=1$, i.e. there is only one linearly independent complex twistor spinor on $(M,c)$. Such examples are easy to generate, as one might for example take a generic Lorentzian metric admitting a parallel spinor as classified in \cite{br} for low dimensions. 
\begin{Proposition} \label{stuff4}
Suppose that the tractor conformal superalgebra of a simply-connected Lorentzian conformal spin manifold $(M^{1,n-1},c)$ satisfies dim $\g_1=1$. Then $\g=\g_0 \oplus \g_1$ is a Lie superalgebra.
\end{Proposition}

\textit{Proof. }
We fix a nontrivial twistor spinor $\psi \in \g_1$ which is unique up to multiplication in $\C^*$ and assume that $\text{ker }\psi = \{0\}$. By Theorem \ref{bg}, this implies that
\begin{align*}
Hol(M,c) \subset \begin{cases} SU \left(1,\frac{n}{2}\right) & (1) \\ O(r) \times SU(1,\frac{n-r}{2}) & (2)\end{cases}
\end{align*}
However, the action of $SU(1,\frac{n}{2})$ on $\Delta^{\C}_{2,n}$ fixes two spinors (cf. \cite{kath}), which excludes $(1)$. In case $(2)$ we have that the representation $\rho$ of $Hol(\overline{\mathcal{Q}}^1_+,\overline{\widetilde{\omega}}^{nc}) \subset Spin^+(2,n)$ on $\Delta^{\C}_{2,n}$ splits into a product of representations $\rho \cong \rho_1 \otimes \rho_2$ on $Spin^+(0,r) \times Spin^+(2,n-r)$.  Furthermore,
\[ \Delta_{2,n}^{\C} \cong \Delta^{\C}_{0,r} \otimes \Delta^{\C}_{2,n-r}, \]
considered as $Spin^+(0,r) \times Spin^+(2,n-r)$-representations. As there exists a $Hol(\overline{\mathcal{Q}}^1_+,\overline{\widetilde{\omega}}^{nc})$-invariant spinor in $\Delta^{\C}_{2,n}$, we conclude (cf. \cite{ldr}) that each of the factors $\rho_1$ and $\rho_2$ admits an invariant spinor. However, $\left(\rho_2 \right)_*$ is the action of a subalgebra of $\mathfrak{su}(1,\frac{n-r}{2})$ on $\Delta^{\C}_{2,n-r}$ which annihilates at least two linearly independent complex spinors. Consequently, the representation $\rho_2$ fixes at least two linearly independent complex spinors and $\rho_1$ fixes at least one nontrivial complex spinor such that $\rho$ fixes at least two linearly independent complex spinors which means that dim $\g_1 > 1$ , in contradiction to our assumption.
Consequently, we have that ker $\psi \neq \{0\}$ for every $\psi \in \g_1$. The second part of the proof of Theorem \ref{hola} then shows that $\g$ is a Lie superalgebra.
$\hfill \Box$

\begin{kor}
If the tractor conformal superalgebra associated to a simply-connected Lorentzian conformal spin manifold $(M,c)$ is no Lie superalgebra, then there exist at least two linearly independent complex twistor spinors on $(M,c)$.
\end{kor}

\subsection*{The tractor conformal superalgebra of flat Minkowski space}
We describe the even part of the conformal algebra of flat Minkowski space $\R^{1,n-1}$ in terms of conformal tractor calculus and discuss extensions to a superalgebra. 
In physics notation (cf. \cite{msl,raj}), the conformal algebra of Minkowski space $\R^{1,n-1}$ with coordinates $x^i$ and the standard flat metric $g_{ij}$ is generated by $P_i, M_{ij}, D$ and $K_i$ - corresponding to translations, rotations, the dilatation and the special orthogonal transformations:
\begin{align*}
P_i &= \partial_i, \\
M_{ij} &=x_i \partial_j - x_j \partial_i, \\
D &= x^i \partial_i, \\
K_i &= 2x_i x^j \partial_j -g(x,x) \partial_i.
\end{align*} 
The Lie brackets can be found in \cite{raj}. As $\R^{1,n-1}$ is conformally flat, all conformal vector fields are automatically normal conformal, and thus the above vector fields generate the algebra $\mathfrak{X}^{nc}(\R^{1,n-1})=\mathfrak{X}^{c}(\R^{1,n-1})$. 
We now consider the following natural isomorphism:
\begin{align}
\tau_0 : \mathfrak{X}^{nc}(\R^{1,n-1}) \stackrel{g}{\cong} Par\left(\Lambda^2_{\mathcal{T}}\left(\R^{1,n-1}\right),\nabla^{nc}\right) \stackrel{\alpha \mapsto \alpha(0)}{\cong} \mathfrak{so}(2,n), \label{idid}
\end{align}
yielding that for flat Minkowski space $\g_0 \cong \mathfrak{so}(2,n)$ on the tractor level. Solving the twistor equation on $\R^{1,n-1}$ is straightforward (cf. \cite{bfkg}): We have for a twistor spinor $\ph \in \Gamma(\R^{1,n-1},S^g_{\C}) \cong C^{\infty}\left(\R^{1,n-1},\Delta_{1,n-1}^{\C}\right)$ using $K^g=0$ that $\nabla D^g \ph = 0$. Consequently, $D^g \ph =: \ph_1 $ is a constant spinor. Integrating the twistor equation along the line $\{ s \cdot x \mid 0 \leq s \leq 1 \}$ yields  that $\ph(x) - \ph(0)= - \frac{1}{n} x \cdot \ph_1$. Thus, $\ph$ is of the form $\ph (x)= \ph_0 - \frac{1}{n} x \cdot \ph_1$. Clearly, this establishes an isomorphism
\begin{align*}
\begin{array}{cccccccc}
\tau_1 : & \text{ker }P^g & \rightarrow & \Delta^{\C}_{1,n-1} \oplus  \Delta^{\C}_{1,n-1} & \cong & {Par}\left( \mathcal{S}(\R^{1,n-1}),\nabla^{nc}\right) & \cong & \Delta_{2,n}^{\C}, \\
& \ph & \mapsto & \left(\ph_0, - \frac{1}{n} \ph_1 \right) & \mapsto & \psi:= \left(\widetilde{\Phi}^g \right)^{-1} \left(\ph_0, - \frac{1}{n} \ph_1 \right) & \mapsto & \psi(0).
\end{array}
\end{align*}
Consequently, the tractor conformal superalgebra of $\R^{1,n-1}$is nothing but $\Lambda^2_{2,n} \oplus \Delta_{2,n}^{\C}$ with brackets as introduced in section \ref{cts}. By means of $\tau_1$ and $\tau_2$ we have an identification
\begin{align}
\g \cong \Lambda^2_{2,n} \oplus \Delta_{2,n}^{\C} \stackrel{\tau_0,\tau_1} \cong \mathfrak{X}^{nc} (\R^{1,n-1}) \oplus \text{ker }P^g, \label{sss}
\end{align}
and the right hand side of (\ref{sss}) is precisely the conformal superalgebra of Minkowski space wrt. the fixed standard metric as considered in \cite{raj,med}, for example, emphasising that the tractor approach to conformal superalgebras is equivalent to the classical approaches. Using an explicit Clifford representation, one directly calculates that $\g$ is no Lie superalgebra if $n > 3$, as also follows from Theorem \ref{hola}. In case $n=3$, and considering Minkowski space $\R^{2,1}$, there is a real structure on $\Delta^{\C}_{3,2}$ (cf. \cite{lm,br}), and restricting ourselves to real twistor spinors leads to the Lie superalgebra\footnote{The odd-odd-odd Jacobi identity holds in this case as every nonzero spinor $v \in \Delta_{3,2}^{\R}$ is pure, from which $\alpha^2_{v,v} \cdot v = 0$ follows. Note that there is no real structure on $\Delta_{2,3}^{\R}$.}
\[ \mathfrak{X}^{nc} (\R^{2,1}) \oplus \text{ker }P_{\R}^g \cong \Lambda^2_{3,2} \oplus \Delta_{3,2}^{\R} \subset  \Lambda^2_{3,2} \oplus \Delta_{3,2}^{\C} = \mathfrak{X}^{nc} (\R^{2,1}) \oplus \text{ker }P_{\C}^g. \]

\section{A tractor conformal superalgebra with R-symmetries for Fefferman spaces} \label{rsymme}
Our aim is to reproduce the construction of conformal Lie superalgebras with R-symmetries for Fefferman spaces as known from \cite{med} in the framework of the conformal tractor calculus. Let $(M,c)$ be an simply-connected, even-dimensional Lorentzian conformal spin manifold and $g \in c$. For the definition and construction of Fefferman spin spaces as total spaces of $S^1-$bundles over strictly pseudoconvex manifolds we refer to \cite{bafe,bl,baju}. The following is a standard fact:

\begin{Proposition}[\cite{lei,bafe}] \label{tuefo}
On a Lorentzian Fefferman spin space $(M^{1,n-1},g)$ there are distinguished, linearly independent complex twistor spinors $\ph_{\epsilon}$ for $\epsilon = \pm 1$ such that 
\begin{enumerate}
\item The Dirac current $V_{\ph_{\epsilon}}$ is a regular lightlike Killing vector field.
\item $\nabla_{V_{\ph_{\epsilon}}} \ph_{\epsilon} = i c \ph_{\epsilon}$ for some $c \in \R \backslash \{0 \}$.
\end{enumerate}
\end{Proposition}
In fact, a Fefferman spin space can also be equivalently described by the existence of twistor spinors $\ph_{\epsilon}$ with Dirac current $V_{\ph_{\epsilon}}$ satisfying 1. and 2. In terms of conformal holonomy, a Fefferman metric in the conformal class is characterized by $Hol(M,c) \subset SU(1,\frac{n}{2})$, cf. \cite{baju,leihabil}.
We now restrict ourselves to \textit{generic} Fefferman spin spaces, i.e. our overall assumption in this section in terms of conformal data is
\begin{align*}
Hol(M^{1,n-1},c) \subset SU\left(1,\frac{n}{2}\right) \text{ and } \text{dim}_{\C}\text{ ker }P^g = 2.
\end{align*}
In case $Hol(M^{1,n-1},c=[g]) = SU\left(1,\frac{n}{2}\right)$, the second requirement follows automatically.

\begin{Proposition} \label{pry}
For Lorentzian conformal structures with $Hol(M,c)=SU\left(1,\frac{n}{2}\right)$ one has that $\text{dim}_{\C}\g_1=2$ and the tractor conformal algebra is \textbf{no} Lie superalgebra.
\end{Proposition}
\textit{Proof. }In order to prove this Proposition, we start with the observation that by (\ref{sytr}) complex parallel spin tractors on $M$ correspond (after fixing a basepoint) to spinors in $\Delta_{2,n}^{\C}$ which are annihilated by the action of $\lambda_*^{-1}\left(\mathfrak{su}\left(1,\frac{n}{2}\right)\right)$. Let us call the space of these spinors $V_{\mathfrak{su}}$. We fix the following \textit{complex} representation of the complex Clifford algebra $Cl^{\C}_{2,n}$ with $n+2=:2m$ on $\C^{2^m}$ (cf. \cite{ba81}):
Let $E,D,U$ and $V$ denote the $2 \times 2$ matrices
\begin{align*} 
E = \begin{pmatrix} 1 & 0 \\ 0 & 1 \end{pmatrix} \text{ , } D = \begin{pmatrix} 0 & -i \\ i & 0 \end{pmatrix} \text{ , } U = \begin{pmatrix} i & 0 \\ 0 & -i \end{pmatrix} \text{ , } V = \begin{pmatrix} 0 & i \\ i & 0 \end{pmatrix}.
\end{align*}
Furthermore, let $\tau_j =\begin{cases}  1 &  \epsilon_j = 1, \\ i & \epsilon_j = -1. \end{cases}$. ${Cl}^{\C}(p,q) \cong M_{2^m}(\C)$ as complex algebras, and an explicit realisation of this isomorphism is given by
\begin{align*}
\Phi_{p,q} (e_{2j-1})&= \tau_{2j-1} \cdot E \otimes...\otimes E \otimes U \otimes \underbrace{D \otimes...\otimes D}_{(j-1) \times},\\
\Phi_{p,q} (e_{2j})  &= \tau_{2j} \cdot E \otimes...\otimes E \otimes V \otimes \underbrace{D \otimes...\otimes D}_{(j-1) \times}.
\end{align*}

We set $\widetilde{u}(\epsilon):=\frac{1}{\sqrt{2}}\cdot \begin{pmatrix} 1 \\ -i \epsilon \end{pmatrix}$ for $\epsilon = \pm 1$ and introduce the spinors
$\widetilde{u}(\epsilon_m,....,\epsilon_1):=\widetilde{u}(\epsilon_m) \otimes...\otimes \widetilde{u}(\epsilon_1)$ which form a basis of $\Delta_{2,n}^{\C}$.
We work with the $Spin^+(2,n)$-invariant scalar product $\langle u, v \rangle_{\Delta^{\C}_{2,n}}:=i (e_1\cdot e_2 \cdot u,v )_{\C^{2^m}}$. One calculates that
\begin{align} \label{625}
\langle \widetilde{u}(\epsilon_m,...,\epsilon_1), \widetilde{u}(\delta_m,...,\delta_1) \rangle = \begin{cases} 0 & (\epsilon_m,...,\epsilon_1) \neq (\delta_m,...,\delta_1) \\
\epsilon_1 & (\epsilon_m,...,\epsilon_1) = (\delta_m,...,\delta_1) \end{cases}
\end{align}
It is now straightforward to compute (cf. \cite{kath}) that 
\begin{align} \label{asint}
V_{\mathfrak{su}}:= \{ v \in \Delta_{2,n}^{\C} \mid \lambda_*^{-1}\left(\mathfrak{su}\left(1,\frac{n}{2}\right)\right) \cdot v = 0 \} = \text{span}_{\C} \{u_+:=\widetilde{u}(1,...,1), u_-:=\widetilde{u}(-1,...,-1) \}.
\end{align}

Another straightforward computation involving (\ref{bla}) and (\ref{625}) yields that \begin{align}
\alpha_{u_{\pm}}^2= \sum_{i=1}^{\frac{n}{2}+1}\epsilon_{2i} \cdot e_{2i-1}^{\flat} \wedge e_{2i}^{\flat}, \label{asfo}
\end{align}
from which follows that 
\begin{align}
\alpha^2_{u_+} \cdot u_+ = i\cdot \left(\frac{n}{2}-1\right) u_+ \neq 0. \label{tgg}
\end{align}
If we turn to geometry, a global version of the previous observations shows that for simply-connected conformal structures with irreducible holonomy $SU(1,\frac{n}{2})$ the dimension of the complex space of twistor spinors is two-dimensional and (\ref{tgg}) yields that the tractor conformal superalgebra is no Lie superalgebra as the odd-odd-odd Jacobi identity is not satisfied.
$\hfill \Box$ 


\subsection*{Algebraic preparation}
We want to investigate the space of parallel spin tractors on $(M,c)$ more closely. To this end, we use the complex spinor representation on $\Delta_{2,n}^{\C}$ from the proof of Proposition \ref{pry} with distinguished spinors $u_{\pm}$. Let $W:=\text{span}_{\C} \{u_+,u_-\}$. We have already computed $\omega_0:= \alpha_{u_{\pm},u_{\pm}}^2$ in (\ref{asfo}). A straightforward, purely algebraic calculation reveals the following:

\begin{Proposition} \label{strata0}
The pseudo-K\"ahler form $\omega_0$ on $\R^{2,n}$ is distinguished by the following properties:
\begin{enumerate}
\item For every $w \in W$ there exists a constant $c_w \geq 0$ such that $\alpha^2_{w,w} = c_w \cdot \omega_0$.
\item $||\omega_0||_{2,n}^2 = \frac{n}{2}+1$
\end{enumerate}
\end{Proposition}
Moreover, one calculates that for all $a,b \in \C$ we have $\frac{1}{i}\omega_0 \cdot (a u_+ + b u_-) = \left(\frac{n}{2}-1\right)\cdot (au_+ - b u_-)$, whence
\begin{align}
\text{span}_{\C} \{ u_{\pm} \} = \text{Eig}_{\C} \left(\frac{1}{i}\omega_0,\pm \left(\frac{n}{2}-1\right)\right). \label{ev}
\end{align}

\begin{Lemma} \label{stuckd}
Consider  $u_+ \in W$ and let $\alpha \in \Lambda^2_{2,n}$ be a 2-form. If $\alpha \cdot u_+ \in W$, then $\alpha$ can be written as $\alpha = \sum_{i=1}^{\frac{n+2}{2}} a_i \cdot e_{2i-1}^{\flat} \wedge e_{2i}^{\flat}$ for $a_i \in \R$. We denote the space of all these forms by $V$.
\end{Lemma}

\textit{Proof. }
We write a generic 2-form as $\alpha = \sum_{i <j} a_{ij} e_{i}^{\flat} \wedge e_j^{\flat}$. It follows that $\alpha \cdot u_{\pm} = \sum_{i<j} a_{ij} e_i \cdot e_j \cdot u_{\pm}$. Using our concrete realisation of Clifford multiplication, one calculates that $j \neq i+1 \Rightarrow e_i e_j u_+ \propto u(1,...,1,-1,1...,1,-1...,1)$, where $-1$ occurs at positions $\left\lfloor  \frac{i+1}{2} \right\rfloor$ and $\left\lfloor  \frac{j+1}{2} \right\rfloor$. As $\alpha \cdot u_+ \in W$, it follows that $a_{ij}=0$ for these choices of $i$ and $j$.
$\hfill \Box$\\

Another purely algebraic computation reveals the following:
\begin{Lemma} \label{strata}
On $W$ there exists a up to sign unique $\C$-linear map $\iota : W \rightarrow W$ such that $\iota^2 = 1$ and $\iota$ is an anti-isometry of $(W,\langle \cdot, \cdot \rangle_{\Delta_{2,n}^{\C}})$, i.e. $\langle \iota(u),\iota(v) \rangle_{\Delta_{2,n}^{\C}} = - \langle u,v \rangle_{\Delta_{2,n}^{\C}}$.
\end{Lemma}
Moreover, (\ref{ev}) shows that setting
\begin{align}
\frac{1}{i}\omega_0 \cdot u =: \left(\frac{n}{2}-1\right) \cdot l(u), \text{ for }u \in W
\end{align}
defines a unique $\C-$linear map $l: W \rightarrow W$. $l$ is an isometry wrt. $\langle \cdot, \cdot \rangle_{\Delta_{2,n}^{\C}}$ and $l^2 = 1$. We note that wrt. the basis $(u_+,u_-)$ of $W$, $\iota$ and $l$ are given by $\iota = \begin{pmatrix} 0 & 1 \\ 1 & 0 \\ \end{pmatrix}$ and $l = \begin{pmatrix} 1 & 0 \\ 0 & -1 \\ \end{pmatrix}$. One easily calculates that for all $\alpha \in V$ and $u \in W$ we have
\begin{equation} \label{calcd}
\begin{aligned}
\alpha \cdot \iota(u) &= - \iota (\alpha \cdot u) \text{, } \alpha \cdot l(u) = l (\alpha \cdot u),\\ \iota(l(u)) &= -l(\iota(u)).
\end{aligned}
\end{equation}

\subsection*{Geometric construction}
We now turn to geometry again: Let $(M,c)$ be a simply-connected Lorentzian conformal spin manifold with special unitary conformal holonomy and suppose that $\text{dim}_{\C} W = 2$, where now $W={Par}(\mathcal{S}_{\C}(M), \nabla^{nc})$. Global versions of our previous algebraic observations show: There exists a unique parallel tractor 2-form $\omega_0 \in {Par}\left(\Lambda^2_{\mathcal{T}}(M),\nabla^{nc} \right)$ distinguished by properties of Proposition \ref{strata0}. Furthermore, Clifford multiplication with $\frac{1}{i}\omega_0$ is an automorphism of $W$ with eigenvalues $\pm (\frac{n}{2}-1)$. We now fix $\psi_{\pm} \in \text{Eig} \left( \frac{1}{i}\omega_0,\pm \left(\frac{n}{2}-1 \right)\right) \cap W$ with $\langle \psi_{\pm},\psi_{\pm} \rangle_{\Delta_{2,n}^{\C}} = \pm 1$ and $\langle \psi_{\pm},\psi_{\mp} \rangle_{\Delta_{2,n}^{\C}} =0$. With these requirements, $\psi_{\pm}$ are unique up to multiplications with elements of $S^1 \subset \C$. In fact, if one fixes a Fefferman metric $g$ in the conformal class, then $\widetilde{\Phi}^g(proj^g_+ \psi_{\pm}) = \ph_{\pm} \in \text{ker }P^g$(up to constant multiples), where $\ph_{\pm}$ were introduced in Proposition \ref{tuefo}. We further require that $\iota(\psi_+)=\psi_-$ which reduces the ambiguity in choosing $\psi_{\pm}$ to only one complex phase. We set 
\begin{align*}
\g_1 :=W=\text{span}_{\C} \{ \psi_+,\psi_-\} \subset \mathcal{S}_{\C}(M).
\end{align*}
On $\g_1$ there are natural maps $\iota : \g_1 \rightarrow \g_1$ and $l: \g_1 \rightarrow \g_1$ with the same properties as the corresponding maps from the algebraic preparations. $\g_1$ defines the odd part of the tractor superconformal algebra we are about to construct. For the construction of the even part, we first set as in section \ref{cts} $\g_0 :=Par \left( \Lambda^2_{\mathcal{T}}(M),\nabla^{nc}\right)$ and equip it with the bracket of endomorphisms.

\begin{Proposition}
For $Hol(M,c) \subset SU\left( 1, \frac{n}{2} \right)$ and dim $\g_1 = 2$, we have that $\g_0$ is abelian and dim $\g_0 \leq \frac{n}{2}+1$.
\end{Proposition}

\textit{Proof. }
For $\alpha \in \g_0$ and $\psi=\psi_+ \in \g_1$ we must by the derivation property of $\nabla^{nc}$ wrt. Clifford multiplication have that also $ \alpha \cdot \psi$ is a parallel spin tractor, i.e. $\alpha \cdot \psi \stackrel{!}{\in} \g_1$. Lemma \ref{stuckd} then determines all possible forms of $\alpha$ and from this the statement is immediate.
$\hfill \Box$\\

We now set $\widetilde{\g}_0:=\g_0 \oplus \R$, where the sum is a direct sum of abelian Lie algebras and thus $\widetilde{\g}_0$ is abelian too. We introduce further brackets on $\g:= \widetilde{\g}_0\oplus \g_1$:
\begin{equation} \label{rbracket}
\begin{aligned}
\left[\cdot, \cdot \right] : \widetilde{\g}_0 \otimes \g_1 & \rightarrow \g_1,\\
\left((\alpha,a),\psi \right) & \mapsto \frac{1}{i} \cdot \left(\alpha \cdot (\iota(\psi))\right) + a \cdot \iota (l(\psi)),\\
[ \cdot, \cdot] : {\g_1} \otimes \g_1 & \rightarrow \widetilde{\g}_0,\\
(\psi_1,\psi_2) & \mapsto \left( \alpha^2_{\psi_1,\psi_2},\left(\frac{n}{2}-1\right) \cdot \text{Re} \left(\langle \psi_1, l(\psi_2) \rangle_{\mathcal{S}} \right) \right).
\end{aligned}
\end{equation}
Clearly, these brackets have the right symmetry properties to turn $\widetilde{\g}_0 \oplus \g_1$ into a superalgebra.
\begin{satz}
The superalgebra $\widetilde{\g}_0 \oplus \g_1$ associated to $(M,c)$ canonically up to sign\footnote{In fact, defining the above brackets via the abstract maps $\iota$ and $l$ rather than using the basis $\psi_{\pm}$ reveals that the construction involves no further choices once $\g_0$ and $\g_1$ are determined.} is a Lie superalgebra. 
\end{satz}
\textit{Proof. }All we have to do is checking the Jacobi identities: By polarization, the odd-odd-odd Jacobi identity is equivalent to $[[\psi,\psi]\psi]=0$ for all $\psi \in \g_1$. By definition, we have for $\psi=a \psi_+ + b \psi_-$ with $a,b \in \C$ that
\begin{align*}
[[\psi,\psi],\psi] &=\left[\left(\alpha^2_{\psi,\psi},\left(\frac{n}{2}-1\right) \cdot \text{Re} \left(\langle \psi, l(\psi) \rangle_{\mathcal{S}}\right) \right),\psi\right] \\
&=\frac{1}{i} \cdot \alpha_{\psi,\psi}^2 \cdot \iota(\psi) + \left(\frac{n}{2}-1\right) \cdot \text{Re} \left(\langle \psi, l(\psi) \rangle_{\mathcal{S}}\right) \cdot \iota (l (\psi))\\
&=\frac{1}{i} \left(|a|^2 \omega_0 + |b|^2 \omega_0 \right) \cdot (a \psi_- + b \psi_+) + \left(\frac{n}{2}-1\right) (|a|^2+|b|^2) \cdot (a \psi_- - b \psi_+) \\
&=(|a|^2 + |b|^2) \cdot \left(\frac{n}{2}-1 \right) \cdot (-a \psi_- + b \psi_+) + \left(\frac{n}{2}-1\right) (|a|^2 + |b|^2) (a \psi_- - b \psi_+) \\
&=0.
\end{align*}
As $\widetilde{\g}_0$ is abelian, the even-odd-odd identity is by polarization equivalent to $\left[\left[(\alpha, \gamma),\psi \right], \psi \right] = 0$ for all $\alpha \in \g_0$, $\gamma \in \R$ and $\psi \in \g_1$. By definition of the brackets involved, this is the case iff
\begin{align} \label{g0part}
\left( \alpha^2_{\frac{1}{i} \cdot \left(\alpha \cdot \iota(\psi) \right) + \gamma \cdot (\iota(l(\psi))), \psi}, \left(\frac{n}{2}-1\right) \cdot \text{Re} \left(\langle \frac{1}{i} \cdot \left(\alpha \cdot \iota(\psi) \right) + \gamma \cdot (\iota(l(\psi))) , l(\psi) \rangle_{\mathcal{S}} \right) \right) \stackrel{!}{=} 0 \in \g_0 \oplus \R.
\end{align}
We again write $\psi = a \psi_+ + b \psi_-$ for complex constants $a$ and $b$. Lemma \ref{stuckd} implies that $\frac{1}{i}\alpha \cdot {\psi}_+ = d \cdot \psi_+$ for some real constant $d$ and $\frac{1}{i} \alpha \cdot \psi_- = -d \cdot \psi_-$. Then the $\g_0-$ part of (\ref{g0part}) is given by
\begin{align*}
\alpha^2_{-ad \psi_- + bd \psi_+ + \gamma(a \psi_- - b \psi_+), a\psi_+ + b \psi_-} = 0,
\end{align*}
where we used that $\alpha^2_{\psi_+,\psi_-}=0$, and the $\R-$part of (\ref{g0part}) is proportional to
\begin{align*}
&  \text{Re}\left(\langle \frac{1}{i} \cdot \left(\alpha \cdot \iota(\psi) \right) + \gamma \cdot (\iota(l(\psi))) , l(\psi) \rangle_{\mathcal{S}}\right) \\
= & \text{Re}\left(\langle d \cdot \left(-a \cdot \psi_- + b \cdot \psi_+ \right) + \gamma \cdot (a \psi_- - b \psi_+), a\psi_+ - b \psi_- \rangle_{\mathcal{S}}\right) \\
=&0.
\end{align*}
Finally, since $\widetilde{\g}_0$ is abelian, the even-even-odd identity is equivalent to
\begin{align*}
\left[(\alpha,a),\frac{1}{i} \beta \cdot (\iota (\psi)) + b \cdot \iota(l(\psi)) \right] \stackrel{!}{=} \left[(\beta,b),\frac{1}{i} \alpha \cdot (\iota (\psi)) + a \cdot \iota(l(\psi)) \right] \in \g_1, 
\end{align*}
where $(\alpha,a), (\beta,b) \in \widetilde{\g}_0$ and $\psi \in \g_1$. Unwinding the definitions and using (\ref{calcd}), we find that the left hand side is given by
\begin{align*}
&\frac{1}{i}\left(\frac{1}{i} \alpha \cdot \iota (\beta \cdot \iota (\psi)) + {a} \cdot \iota (l(\beta \cdot \iota (\psi))) + b \cdot \alpha \cdot l(\psi) \right)+ ab \cdot \iota(l(\iota(l(\psi)))) \\
{=}& \alpha \cdot \beta \cdot \psi + \frac{a}{i} \beta \cdot l(\psi) + \frac{b}{i} \cdot \alpha \cdot l(\psi) - ab \cdot \psi \stackrel{[\alpha, \beta]=0}{=} \beta \cdot \alpha \cdot \psi + \frac{a}{i} \beta \cdot l(\psi) + \frac{b}{i} \cdot \alpha \cdot l(\psi) - ab \cdot \psi \\
=& \left[(\beta,b),\frac{1}{i} \alpha \cdot (\iota (\psi)) + a \cdot \iota(l(\psi)) \right].
\end{align*}
These calculations prove the Theorem.
$\hfill \Box$

\begin{bemerkung}
Let $g \in c$ be a Fefferman metric on $M$. By means of $g$ we identify the parallel spin tractors $\psi_{\epsilon}$ with the distinguished twistor spinors $\ph_{\epsilon}$ from Proposition \ref{tuefo} for $\epsilon = \pm 1$ and parallel 2-form tractors with normal conformal vector fields. Calculations completely analogous to that in section \ref{dem} reveal that the even-odd bracket (\ref{rbracket}) is under this $g$-metric identification given by (extension of)
\begin{align*}
(\mathfrak{X}^{nc}(M) \oplus \R) \times \text{ker }P^g & \rightarrow \text{ker }P^g, \\
((V,a),\ph_{\epsilon}) & \mapsto L_{V} \ph_{-\epsilon} + \epsilon \cdot a \cdot \ph_{-\epsilon},
\end{align*}
and in this picture the odd-odd-odd Jacobi identity for $\widetilde{\g}_0 \oplus \g_1$ is equivalent to the existence of a constant $\rho$ such that  $L_{V_{\ph_\epsilon}} \ph_{\epsilon} + \epsilon \cdot \rho \cdot \ph_{\epsilon} = 0$, as proved independently in \cite{med}.
\end{bemerkung}

\begin{bemerkung} \label{oda} There is an odd-dimensional analogue of this construction: Namely, consider the case of a simply-connected, Lorentzian Einstein-Sasaki manifold $(M^{1,n-1},g)$ of negative scalar curvature (cf. \cite{lei,boh}), which can be equivalently characterized in terms of special unitary holonomy of the cone over $(M,g)$. It follows that $(M,g)$ is spin and there again exist two distinguished conformal Killing spinors (cf. \cite{bl,leihabil}). Let us assume that the complex span of these twistor spinors is already $\text{ker }P^g=: \g_1$. As $(M,g)$ is Einstein with $\text{scal}^g <0$ there exists in this case (cf. \cite{leihabil, baju}) a \textit{distinguished} spacelike, parallel standard tractor $\tau$, defining a holonomy reduction $Hol(M,[g]) \subset SU\left(1,\frac{n-1}{2}\right) \subset SO(2,n-1) \subset SO(2,n)$ and a splitting $\mathcal{T}(M) = \langle \tau \rangle^{\bot} \oplus \langle \tau \rangle$. Furthermore $\Delta_{2,n-1} \cong \Delta_{2,n}$ as $Spin(2,n-1)$-representations. Setting $\g_0 := {Par}\left(\Lambda^2_{\mathcal{T}}(M),\nabla^{nc} \right) \cap \{ \alpha \in \Omega^2_{\mathcal{T}}(M) \mid \alpha(\tau,\cdot) = 0 \}$, we can then proceed completely analogous to the even-dimensional case just discussed, i.e. we perform in the tractor setting the same purely algebraic construction on the orthogonal complement of $\tau$ in $\mathcal{T}(M)$. This turns $(\g_0 \oplus \R) \oplus \g_1$ into a Lie superalgebra with R-symmetries. Again, the overall construction is canonical. For a construction which uses a fixed metric in the conformal class, we refer to \cite{med}.
\end{bemerkung}

\section{Summary and application in small dimensions} \label{66}
We want to summarize the various possibilities and obstructions one faces in the attempt of constructing a conformal Lie superalgebra via the tractor approach in Lorentzian signature. To this end, recall that twistor spinors on Lorentzian manifolds can be categorized into three types according to Theorem \ref{bg}: We have shown:
\begin{itemize}
\item If all twistor spinors are of type 1. or 2., the tractor conformal superalgebra is a Lie superalgebra (cf. Theorem \ref{hola}). Moreover, if $\g$ is a Lie superalgebra, there is a Brinkmann metric in the conformal class or a local splitting $[g]=[-dt^2+h]$, where $h$ is Riemannian Ricci-flat K\"ahler.
\item The previous situation always occurs if the space of twistor spinors is 1-dimensional.
\item If there are exactly two linearly independent twistor spinors of type $3.a$ or $3.b$ (depending on the dimension to be even or odd), one can construct a Lie superalgebra under the inclusion of an R-symmetry. Depending on the dimension, one has a Fefferman metric or a Lorentzian Einstein-Sasaki metric in the conformal class.
\end{itemize}

\begin{bemerkung}
We have not yet discussed the case when the twistor spinor is of type $3.c$ in Theorem \ref{bg}, i.e. when there is -at least locally - a splitting $(M,g) \cong (M_1,g_1) \times (M_2,g_2)$ into a product of Einstein spaces. By \cite{leihabil,baju} we have that $Hol(M,[g]) \cong Hol(M_1,[g_1]) \times Hol(M_2,[g_2])$. In this situation, it is an algebraic fact (cf. \cite{ldr}) that every spinor $v \in \Delta_{2,n}$ which is fixed by $Hol(M,[g])$ is of the form $v=v_1 \otimes v_2$ where $Hol(M_i,[g_i])v_i = v_i$. As also the converse is trivially true, we see that on the level of tractor conformal superalgebras, the product case manifests itself in a splitting of the odd part of $\mathfrak{g}$, i.e. $\g_1 = \g_1^1 \otimes \g_1^2$, where $\g_1^i$ are the odd parts of the tractor conformal Lie superalgebras $\g^i = \g_0^i \oplus \g_1^i$ of $(M_i,[g_i])$ for $i=1,2$. Moreover, note that we never have a splitting in the even part, $\g_0 \neq \g_0^1 \oplus \g_0^2$. This is because as $(M_i,g_i)$ are Einstein manifolds, there are parallel standard tractors $t_i \in \mathcal{T}(M_i)$, and it follows that $t_1 \wedge t_2 \in \mathfrak{g}_0$, but $t_1 \wedge t_2 \notin \g_0 \oplus \g_1$. It is moreover clear from the structure of $\alpha_{\psi}^2$ from Theorem \ref{bg} in this situation that $\g$ is no Lie superalgebra in this case. \cite{med} presents a way of extending $\g$ to a Lie superalgebra under the inclusion of R-symmetries.
\end{bemerkung}

We have now studied the construction of a tractor conformal superalgebra for every (local) geometry admitting twistor spinors and summarize our results:
\begin{satz} \label{satan}
Let $(M^{1,n-1},c)$ be a Lorentzian conformal spin manifold admitting twistor spinors. Assume further that all twistor spinors on $(M,c)$ are of the same type according to Theorem \ref{bg}. Then there are  the following relations between special Lorentzian geometries in the conformal class $c$ and properties of the tractor conformal superalgebra $\g = \g_0 \oplus \g_1$ of $(M,c)$:
\begin{center}
\begin{tabular}[h]{|p{3cm}|p{4.5cm}|p{5.5cm}|}
  \hline
 Twistor spinor type (Thm. \ref{bg})  & Special geometry in c & Structure of $\g = \g_0 \oplus \g_1$ \\ \hline\hline
1.  & Brinkmann space & Lie superalgebra	 \\ \hline
2.  & Splitting $(\R,-dt^2) \times $ Riem. Ricci-flat & Lie superalgebra \\ \hline
3.a & Lorentzian Einstein Sasaki (n odd) & No Lie superalgebra, becomes Lie superalgebra under inclusion of nontrivial R-symmetry \\ \hline
3.b & Fefferman space (n even) & No Lie superalgebra, becomes Lie superalgebra under inclusion of nontrivial R-symmetry \\ \hline
3.c & Splitting $M_1 \times M_2$ into Einstein spaces & No Lie superalgebra, odd part splits $\g^i = \g_0^i \otimes \g_1^i$, but $\g_0 \neq \g_0^1 \oplus \g_0^2$ \\ \hline
\end{tabular} 
\end{center}
\end{satz}

Let us apply this statement to tractor conformal superalgebras $\g$ of non conformally flat Lorentzian conformal manifolds $(M^{1,n-1},[g])$ admitting twistor spinors in small dimensions which have been studied in \cite{lei, bl}:\\
\newline
\textit{Let n=3. }It is known that dim ker $P^g \leq 1$ in this situation. Consequently, by Proposition \ref{stuff4} $\g$ is a tractor conformal Lie superalgebra.  Every twistor spinor is off a singular set locally equivalent to a parallel spinor on a $pp$-wave.\\
\textit{Let n=4. }Here, dim ker $P^g \leq 2$. Exactly one of the following cases occurs: Either, there is a Fefferman metric in the conformal class with two linearly independent twistor spinors. In this case we can construct a tractor superalgebra with R-symmetries. Otherwise, all twistor spinors are locally equivalent to parallel spinors on pp-waves. In this case the ordinary construction of a tractor conformal \textit{Lie} superalgebra $\g$ works.\\
\textit{Let n=5. }This case is already more involved but the possibility of constructing a tractor conformal Lie superalgebra can be completely described: One again has that 
dim ker $P^g \leq 2$. Exactly one of the following cases occurs:
\begin{enumerate}
\item There is a Lorentzian Einstein Sasaki metric in the conformal class. In this case, dim ker $P^g = 2$ and one can construct a tractor conformal Lie superalgebra with R-symmetries as indicated in Remark \ref{oda}.
\item $(M,g)$ is (at least locally) conformally equivalent to a product $\R^{1,0} \times (N^4,h)$, where the last factor is Riemannian Ricci-flat K\"ahler and admits two linearly independent parallel spinors. This corresponds to type 2. twistor spinors from Theorem \ref{bg}, and thus one can construct a tractor conformal Lie superalgebra.
\item All twistor spinors are equivalent to parallel spinors on $pp$-waves. Again, the construction yields a Lie superalgebra.
\end{enumerate}
\textit{Let }$n \geq 6$. Now \textit{mixtures} can occur, i.e. it is possible that some twistor spinors are of type 1. or 2., and some twistor spinors are of type 3. according to Theorem \ref{bg}. In this case, Theorem \ref{satan} does not apply.

\section{Extension to higher signatures} \label{highersgn}
Let $(M^{p,q},c)$ be a space- and time oriented conformal spin manifold of arbitrary signature $(p,q)$ with $p+q=n$ and complex space of parallel spin tractors $\g_1 = Par(\mathcal{S}_{\C}(M),\nabla^{nc})$. We want to associate to $(M,c)$ a tractor conformal superalgebra in a natural way. However, our construction from the previous sections depends crucially on Lorentzian signature. More precisely, the bracket $\g_1 \times \g_1 \rightarrow \g_0$ may become trivial in other signatures, and it therefore has to be modified: Every parallel spin tractor on $M$ naturally gives rise to a series of parallel tractor $k-$forms, which are nontrivial at least for $k=p+1$ . We include all these conformal symmetries in the algebra and use them to construct the odd-odd-bracket. Thus we then also have to modify $\g_0$ and would like to set
\begin{align}
\g_0 := Par\left(\Lambda^*_{\mathcal{T}}(M),\nabla^{nc} \right) \subset \Omega^*_{\mathcal{T}}(M). \label{atartup}
\end{align}

\subsection*{Algebraic preparation}
Let us for a moment change our notation to $\R^{r,s}$ and $m=r+s$, as the following results will later be applied to $\R^{p,q}$ and $\R^{p+1,q+1}$. In order to introduce a bracket on $\Lambda^k_{r,s}$, we recall the following formulas for the action of a vector $X \in \R^{r,s}$ and a $k$-form $\omega \in \Lambda^k_{r,s}$ on a spinor $\ph \in \Delta^{\C}_{r,s}$ (cf. \cite{bfkg}):
\begin{equation} \label{ext1}
\begin{aligned}
X \cdot (\omega \cdot \ph) &= (X^{\flat} \wedge \omega) \cdot \ph - (X \invneg \omega) \cdot \ph, \\
\omega \cdot (X \cdot \ph) &= (-1)^k \left((X^{\flat} \wedge \omega) \cdot \ph + (X \invneg \omega) \cdot \ph \right). 
\end{aligned}
\end{equation}
This motivates us to set $X \cdot \omega := X^{\flat} \wedge \omega - X \invneg \omega \in \Lambda^{k-1} \oplus \Lambda^{k+1}$ for $X \in \R^{r,s}$ and $\omega \in \Lambda^k_{r,s}$. We use this to set inductively for $e_{I}^{\flat}:=e^{\flat}_{i_1} \wedge ...\wedge e^{\flat}_{i_j} \in \Lambda^j_{r,s}$, where $1 \leq i_1 < i_2<...<i_j \leq n$:
\begin{align}
e_{I}^{\flat} \cdot \omega := e_{i_1} \cdot (e^{\flat}_{I \backslash \{i_1 \}} \cdot \omega ). \label{chapp}
\end{align}
By multilinear extension, this defines $\eta \cdot \omega \in \Lambda^*_{r,s}$ for all $\eta, \omega \in \Lambda^*_{r,s}$. One checks that this product is associative and $O(r,s)-$equivariant, i.e.
\begin{align} \label{ete}
(A\eta) \cdot (A\omega) = A(\eta \cdot \omega) \text{ }\forall A \in O(r,s).
\end{align}
\begin{bemerkung}
The above definition of $\cdot$ is useful for concrete calculations. However, there is an equivalent way of introducing the inner product $\cdot$ on the space of forms which shows that this construction is very natural. To this end, consider the multilinear maps
\[ f_k: \underbrace{\R^{r,s} \times...\times \R^{r,s}}_{k \text{ times}} \rightarrow Cl_{r,s}, \text{ }(v_1,...v_k) \mapsto \frac{1}{k!} \sum_{\sigma \in S_k} \text{sgn}(\sigma) \cdot v_{\sigma_1}\cdot...\cdot v_{\sigma_k}. \]
The maps $f_k$ induce a canonical vector space isomorphism (cf. \cite{lm})
\begin{align*}
\widetilde{f}: \Lambda^*_{r,s} \rightarrow Cl_{r,s},
\end{align*}
for which $\widetilde{f}(v_1^{\flat} \wedge...\wedge v_k^{\flat}) = f_k(v_1,...,v_k)$ holds. It is now straightforward to calculate that our inner product (\ref{chapp}) on $\Lambda^*_{r,s}$ is just the algebra structure which makes $\widetilde{f}$ become an algebra isomorphism, i.e. one has for $\eta, \omega \in \Lambda^*_{r,s}$ that
\begin{align}
\eta \cdot \omega = \widetilde{f}^{-1}\left(\widetilde{f}(\eta) \cdot \widetilde{f}(\omega)\right). \label{prodo}
\end{align}
\end{bemerkung}

With these definitions, the space $\Lambda^*_{r,s}$ together with the map
\begin{align}
[\cdot, \cdot]_{\Lambda}:\Lambda^*_{r,s} \otimes \Lambda^*_{r,s} \rightarrow \Lambda^*_{r,s}\text{, }[\eta,\omega]_{\Lambda}:=\eta \cdot \omega - \omega \cdot \eta \label{dota}
\end{align}
becomes a Lie algebra in a natural way due to associativity of Clifford multiplication.

\begin{bemerkung}
We index this bracket with the subscript $\Lambda$ because on 2-forms there are now the bracket $[\cdot, \cdot]_{\Lambda}$ and the endomorphism-bracket $[\cdot, \cdot]_{\mathfrak{so}}$ from  the previous sections. However, it is straightforward to calculate that $[\cdot, \cdot]_{\Lambda}=2\cdot [\cdot,\cdot]_{\mathfrak{so}}$. Whence these two Lie algebra structures are equivalent. Note that $[\Lambda^k_{r,s},\Lambda^l_{r,s}]$ is in general of mixed degree for $k,l \neq 2$.
\end{bemerkung}

\subsection*{Conformally invariant definition of $\mathfrak{g}$}
Turning to geometry again, let $(M^{p,q},c)$ be a conformal spin manifold of signature $(p,q)$. Given $\alpha, \beta \in \Omega^*_{\mathcal{T}}(M)$ and $x \in M$, we may write $\alpha(x)=[s,\widehat{\alpha}]$ and $\beta(x)=[s,\widehat{\beta}]$ for some $s \in \overline{\mathcal{P}}^1_x$ and $\widehat{\alpha}, \widehat{\beta} \in \Lambda^*_{p+1,q+1}$. We then introduce a bracket on tractor forms by setting
\begin{align} \label{pwd}
\left(\alpha \cdot \beta \right)(x):=[s,\widehat{\alpha} \cdot \widehat{\beta}]. 
\end{align}
The equivariance property (\ref{ete}) shows that (\ref{pwd}) is well-defined. We furthermore define the bracket $[\alpha,\beta]_{\mathcal{T}}$ on $ \Omega^*_{\mathcal{T}}(M)$ by pointwise application of (\ref{dota}). Clearly, this defines a Lie algebra structure on $\Omega^*_{\mathcal{T}}(M)$.

\begin{Lemma}
The normal conformal Cartan connection $\nabla^{nc}$ on $\Omega^*_{\mathcal{T}}(M)$ is a derivation wrt. the product $\cdot$, i.e.
\[ \nabla^{nc}_X (\alpha \cdot \beta) = \left(\nabla^{nc}_X \alpha\right)\cdot \beta  + \alpha \cdot \left(\nabla^{nc}_X \beta\right) \text{ }\forall \alpha, \beta \in \Omega^*_{\mathcal{T}}(M) \text{ and }X \in \mathfrak{X}(M). \]
\end{Lemma}

\textit{Proof. }
Suppose first that $\alpha = Y^{\flat}$ for some standard tractor $Y \in \Gamma(\mathcal{T}(M))$. We calculate:
\begin{align*}
\nabla_{X}^{nc}(\alpha \cdot \beta) &= \nabla_X^{nc} (Y^{\flat} \wedge \beta - Y \invneg \beta )\\
&=\left(\nabla_X^{nc} Y\right)^{\flat} \wedge \beta + Y^{\flat} \wedge \left(\nabla_X^{nc} \beta \right) - \left(\nabla_X^{nc} Y\right) \invneg \beta - Y \invneg \left(\nabla_X^{nc} \beta \right) \\
&= \left(\nabla_X^{nc} \alpha \right) \cdot \beta + \alpha \cdot \left(\nabla_X^{nc} \beta \right).
\end{align*}
As a next step, let $\alpha \in \Omega^*_{\mathcal{T}}(M)$ be arbitrary. We fix $x \in M$ and a local pseudo-orthonormal frame $(s_0,...,s_{n+1})$ (wrt $\langle \cdot , \cdot \rangle_{\mathcal{T}}$ ) on $\mathcal{T}(M)$ around $x$ such that  $\nabla^{nc} s_i = 0$ for $i=0,...,n+1$ at $x$. Wrt. this frame we write $\alpha = \sum_{I} \alpha_I s_I^{\flat}$ locally around $x$ for smooth functions $\alpha_I$. We apply the above result inductively for $Y = s_i$ to obtain at $x$
\begin{align*}
\nabla^{nc}_X(\alpha \cdot \beta)& = \sum_I \nabla_X^{nc}\left(\alpha_I \cdot s_I \cdot \beta \right) = \sum_I X(\alpha_I) \cdot s_I \cdot \beta + \sum_I \alpha_I s_I \cdot \nabla^{nc}_X \beta \\
&=\left(\nabla_X^{nc}\alpha \right) \cdot \beta + \alpha \cdot \left(\nabla_X^{nc}\beta \right),
\end{align*}
which gives the desired formula. 
$\hfill \Box$

\begin{kor}
$\alpha, \beta \in Par\left(\Lambda^*_{\mathcal{T}}(M),\nabla^{nc}\right)$ implies that $[\alpha, \beta]_{\mathcal{T}} \in Par\left(\Lambda^*_{\mathcal{T}}(M),\nabla^{nc}\right)$. Thus the space $Par\left(\Lambda^*_{\mathcal{T}}(M),\nabla^{nc}\right)$ together with the bracket induced by $[\cdot, \cdot]_{\mathcal{T}}$ is a Lie subalgebra of $\left(\Omega_{\mathcal{T}}(M)^*,[\cdot, \cdot]_{\mathcal{T}}\right)$.
\end{kor}

As a next step we extend the Lie algebra $\g_0$ (cf. (\ref{atartup})) of (higher order) conformal symmetries 
together with the bracket defined above to a tractor conformal superalgebra $\g=\g_0 \oplus \g_1$ in a natural way by setting as before
\[ \g_1 := Par(\mathcal{S}(M),\nabla^{nc}), \]
and introducing the brackets
\begin{equation} \label{brain}
\begin{aligned}
\g_0 \times \g_0 \rightarrow \g_0\text{,  } & (\alpha, \beta) \mapsto [\alpha,\beta]_{\mathcal{T}}, \\
\g_0 \times \g_1 \rightarrow \g_1\text{,  } & (\alpha, \psi) \mapsto \alpha \cdot \psi, \\
\g_1 \times \g_0 \rightarrow \g_1\text{,  } & (\psi, \alpha) \mapsto - \alpha \cdot \psi, \\
\g_1 \times \g_1 \rightarrow \g_0\text{,  } & (\psi_1, \psi_2) \mapsto \sum_{l \in L_p} \alpha^l_{\psi_1, \psi_2}.
\end{aligned}
\end{equation}
Here, $L_p:=\{ l \in \mathbb{N} \mid \psi \mapsto \alpha^l_{\psi,\psi}\text{ not identically }0,\text{ }\alpha^l_{\psi_1, \psi_2} \text{ symm. in } \psi_1,\psi_2\}$, and for given $p$, one always has $p+1 \in L_p$, and thus the brackets have the right symmetry properties.

\begin{bemerkung}
If $p=1$ and we allow only $l=2$ in the last bracket, we recover a tractor conformal superalgebra which is naturally isomorphic to the one constructed in the previous chapter. Thus, the above construction may be viewed as a reasonable extension to arbitrary signatures.
\end{bemerkung}

It is of course natural to ask, as done in the Lorentzian case, under which conditions the tractor conformal superalgebra actually is a \textit{Lie} superalgebra, i.e. we have to check the Jacobi identities: 
\begin{itemize}
\item As $\g_0$ is a Lie algebra, the even-even-even identity is trivial.
\item It holds by construction of the bracket $[\cdot, \cdot]_{\mathcal{T}}$ as extension of (\ref{ext1}) that 
\[[\alpha, \beta]_{\mathcal{T}} \cdot \psi =\alpha \cdot \beta \cdot \psi - \beta \cdot \alpha \cdot \psi. \]
But this is precisely the even-even-odd Jacobi identity.

\item The Jacobi identity for the odd-odd-odd component again leads to \[ \alpha_{\psi,\psi}^l \cdot \psi \stackrel{!}{=}0.\]However, there is no known way of expressing this condition in terms of $Hol(M,c)$ due to the fact that a classification of possible parallel tractor forms induced by twistor spinors is only available for the Lorentzian case.
\item The even-odd-odd Jacobi identity is by polarization equivalent to \begin{align}
[\alpha,\alpha^l_{\psi,\psi}]_{\mathcal{T}}\stackrel{!}{=}2 \cdot \alpha^l_{\alpha \cdot \psi, \psi} \text{ for }l\in L_p. \label{gt}
\end{align} However, this identity fails to hold in general. From an algebraic point of view this is due to the fact that $[\Lambda^k_{p+1,q+1},\Lambda^k_{p+1,q+1}]_{\Lambda} \subset \Lambda^k_{p+1,q+1}$ only if $k=2$. This was precisely the situation we had in the Lorentzian setting. For other values of $k$ and $p$ the definition of $[ \cdot, \cdot ]_{\mathcal{T}}$ leads to additional terms on the left hand side of (\ref{gt}).
\end{itemize}

\subsection*{Example: Generic twistor spinors in signature (3,2)}
Consider a conformal spin manifold $(M,c)$ in signature $(3,2)$ admitting a generic real twistor spinor, cf. \cite{hs1}. This means that there exists a twistor spinor $\ph \in \Gamma(S^g_{\R})$ such that $\langle \ph, D^g \ph \rangle_{S^g} \neq 0$ (which is independent of $g \in c$). Under further generic assumptions on the conformal structure, one has that $Hol(M,c)=G_{2,2} \subset SO^+(4,3)$\footnote{The existence of a generic real twistor spinor always implies $Hol(M,c)\subset G_{2,2}$, cf. \cite{hs1}. The exact conditions for full holonomy $G_{2,2}$ are given in \cite{leinur} in terms of an explicit ambient metric construction whose metric holonomy coincides with $Hol(M,c)$}, where $G_{2,2}$ can also be defined as the stabilizer of a generic 3-form $\omega_0 \in \Lambda^3_{4,3}$ under the $SO^+(4,3)$-action, see \cite{kath}.\\
Under these conditions, there is up to constant multiples exactly one linearly independent real pure spin tractor $\psi \in \Gamma(M,\mathcal{S}_{\R}(M))$, additionally satisfying dim ker $\psi = 0$. All parallel tractor forms on $(M,c)$ are given by the span of $\alpha_{\psi}^3 \in \Omega^3_{\mathcal{T}}(M)$, being pointwise of type $\omega_0$ and $\ast \alpha_{\psi}^3 \in \Omega^4_{\mathcal{T}}(M)$, being pointwise of type $\ast \omega_0$. Thus, the tractor conformal superalgebra of $(M,c)$ is given by
\begin{align*}
\g=\g_0 \oplus \g_1 = \text{span}\{\alpha_{\psi}^3, \ast \alpha_{\psi}^3 \} \oplus \text{span}\{ \psi \}.
\end{align*}
Pure linear algebra in $\R^{4,3}$ reveals that\footnote{See also \cite{kath} for explicit formulas of $\omega_0$ and pointwise orbit representatives of $\psi$.}
\begin{align*}
\alpha_{\psi}^3  \cdot (\ast \alpha_{\psi}^3) &= (\ast \alpha_{\psi}^3) \cdot \alpha_{\psi}^3, \\
\alpha_{\psi}^3 \cdot \psi &= \text{const.}_1 \cdot \psi, \\
(\ast \alpha_{\psi}^3) \cdot \psi &= \text{const.}_2 \cdot \psi,
\end{align*}
where the $\psi-$dependent constants are proportional to $\langle \psi, \psi \rangle_{\mathcal{S}}$ and zero iff $\psi = 0$. These observations directly translate into the following properties of the superalgebra $\g$ with brackets as introduced in (\ref{brain}).
\begin{Proposition}
The tractor conformal superalgebra $\g$ associated to a conformal spin manifold $(M,c)$ in signature $(3,2)$ with $Hol(M,c)=G_{2,2}$ does not satisfy the odd-odd-odd and the even-odd-odd Jacobi identities. Moreover, the even part $\g_0$ is abelian.
\end{Proposition}  

The example underlines that in contrast to the Lorentzian case, tractor conformal superalgebras need not satisfy at least 3 of the 4 Jacobi identities.

\subsection*{Metric description}
As done in the Lorentzian case, we want to compute the brackets of a general tractor conformal superalgebra $\g=\g_0 \oplus \g_1$ wrt. a metric in the conformal class. To this end, let $\alpha \in \Lambda^{k+1}_{p+1,q+1}, \beta \in \Lambda^{l+1}_{p+1,q+1}$. As in (\ref{mind}) we decompose $\alpha = e_+^{\flat} \wedge \alpha_+ + \alpha_0 + e^{\flat}_- \wedge e_+^{\flat} \wedge \alpha_{\mp} + e_-^{\flat} \wedge \alpha_-$. We want to compute $\left([\alpha,\beta]_{\Lambda}\right)_+$, i.e. the $+$-component of $[\alpha,\beta]_{\Lambda}$ wrt. the decomposition (\ref{mind}). As a preparation, we calculate for $\omega \in \Lambda^{r}_{p,q}, \eta \in \Lambda^s_{p,q}$ the products
{\allowdisplaybreaks \begin{align*}
(e^{\flat}_{\pm} \wedge \omega) \cdot \eta &= e^{\flat}_{\pm} \wedge (\omega \cdot \eta), \\
(e^{\flat}_{\pm} \wedge \omega) \cdot (e_{\pm}^{\flat} \wedge \eta) &= 0, \\
(e^{\flat}_{\pm} \wedge \omega) \cdot (e_{\mp}^{\flat} \wedge \eta) &=(-1)^r \left(e_{\pm}^{\flat} \wedge e_{\mp}^{\flat} \wedge (\omega \cdot \eta) - \eta \cdot \omega \right),\\
(e^{\flat}_{\pm} \wedge \omega) \cdot (e_-^{\flat} \wedge e_+^{\flat} \wedge \eta) &= \mp e_{\pm}^{\flat} \wedge (\omega \cdot \eta), \\
\omega \cdot (e_{\pm}^{\flat} \wedge \eta) &= (-1)^r e_{\pm}^{\flat} \wedge (\omega \cdot \eta), \\
\omega \cdot (e_-^{\flat} \wedge e_+^{\flat} \wedge \eta) &=e_-^{\flat} \wedge e_+^{\flat} \wedge (\omega \cdot \eta), \\
(e_-^{\flat} \wedge e_+^{\flat} \wedge \omega) \cdot \eta &= e_-^{\flat} \wedge e_+^{\flat} \wedge (\omega \cdot \eta), \\
(e_-^{\flat} \wedge e_+^{\flat} \wedge \omega)  \cdot (e_{\pm}^{\flat} \wedge \eta) &= \pm (-1)^r e_{\pm}^{\flat} \wedge (\omega \cdot \eta),\\
(e_-^{\flat} \wedge e_+^{\flat} \wedge \omega) \cdot (e_-^{\flat} \wedge e_+^{\flat} \wedge \eta) &= \omega \cdot \eta.
\end{align*}}
With these formulas, it is straightforward to compute that for $\alpha, \beta$ as above one has
\begin{align}
\left(\alpha \cdot \beta \right)_+ = \alpha_+ \cdot \beta_0 - \alpha_+ \beta_{\mp} +(-1)^{k+1} \alpha_0 \beta_+ + (-1)^{k+1} \alpha_{\mp} \cdot \beta_+, \label{76}
\end{align}
and therefore,
\begin{align}
\left([\alpha,\beta]_{\Lambda}\right)_+ =& \alpha_+ \cdot \beta_0 - (-1)^{l+1} \beta_0 \cdot \alpha_+ - {\alpha_+} \cdot \beta_{\mp} - (-1)^{l+1}\beta_{\mp} \cdot \alpha_+ +(-1)^{k+1} \alpha_0 \cdot \beta_+ - \beta_+ \cdot \alpha_0 \notag \\
&+ (-1)^{k+1} \alpha_{\mp} \cdot \beta_+ + \beta_+ \cdot \alpha_{\mp}. \label{77}
\end{align}
This directly leads to the following global version:

\begin{Proposition}\label{cne}
Let $g \in c$ and let $\alpha, \beta \in \g_0$ be of degree $k+1$ and $l+1$ respectively. Further, let $\alpha_+ = {proj}^g_{\Lambda,+} \alpha \in \Omega^k_{nc,g}(M)$ and  $\beta_+ = {proj}^g_{\Lambda,+} \beta \in \Omega^l_{nc,g}(M)$ denote the associated nc-Killing forms. As $\alpha \cdot \beta \in \g_0$ and $[\alpha,\beta]_{\mathcal{T}} \in \g_0$ are again parallel, the forms  $\left(\alpha \cdot \beta\right)_+= {proj}^g_{\Lambda,+} (\alpha \cdot \beta)\in \Omega^*_{nc,g}(M)$  and  $\left([\alpha,\beta]_{\mathcal{T}}\right)_+ ={proj}^g_{\Lambda,+} ([\alpha,\beta]_{\mathcal{T}})\in \Omega^*_{nc,g}(M)$  are again nc-Killing forms wrt. $g$. They are explicitly given by
\begin{equation}\label{fo1}
\begin{aligned}
\alpha_+ \circ \beta_+ := \left(\alpha \cdot \beta\right)_+ =& \frac{1}{l+1}\cdot \alpha_+ \cdot d\beta_+ + \frac{1}{n-l+1} \alpha_+ \cdot d^* \beta_+ \\
&+ (-1)^{k+1} \frac{1}{k+1} d\alpha_+ \cdot \beta_+ + (-1)^k \cdot \frac{1}{n-k+1} d^* \alpha_+ \cdot \beta_+,
\end{aligned}
\end{equation}
\begin{align*}
\left([\alpha,\beta]_{\mathcal{T}}\right)_+ =& \frac{1}{l+1}\cdot \alpha_+ \cdot d\beta_+ + (-1)^l \frac{1}{l+1} d\beta_+ \cdot \alpha_+ + \frac{1}{n-l+1} \alpha_+ \cdot d^* \beta_+ +(-1)^{l+1} \frac{1}{n-l+1} d^*\beta_+ \cdot \alpha_+   \\ 
&+(-1)^{k+1} \frac{1}{k+1} d\alpha_+ \cdot \beta_+ - \frac{1}{k+1} \beta_+ \cdot d\alpha_+ + (-1)^k \cdot \frac{1}{n-k+1} d^* \alpha_+ \cdot \beta_+ - \frac{1}{n-k+1} \beta_+ \cdot d^* \alpha_+. 
\end{align*}
\end{Proposition}

\textit{Proof. }
This follows directly from the explicit form of the isomorphism ${proj}^g_{\Lambda,+}$ from (\ref{soga}), i.e. one has to insert $(\alpha_+,\alpha_0,\alpha_{\mp},\alpha_-)=(\alpha_+, \frac{1}{k+1}d\alpha_+, -\frac{1}{n-k+1}d^* \alpha_+, \Box_k \alpha_+)$  into the formulas (\ref{76}), (\ref{77}).
$\hfill \Box$\\

We study some interesting consequences and applications. First, note that Proposition (\ref{cne}) opens a way to construct new nc-Killing forms out of existing ones, i.e. $\circ$ defines a map  
\begin{align}
\circ: \Omega^k_{nc,g}(M) \times \Omega^l_{nc,g}(M) \rightarrow \Omega^*_{nc,g}(M). \label{map}
\end{align}
In general, the resulting product is of mixed degree. We have already shown in Proposition \ref{cof} that for nc-Killing 1-forms the bracket $[ \cdot, \cdot ]_{\mathcal{T}}$ corresponds via fixed $g \in c$ to the Lie bracket of vector fields (up to a factor). For deg $\alpha = 2$ one can simplify the expression from Proposition \ref{cne} as follows:

\begin{Proposition} \label{na}
Let $\alpha \in Par(\Lambda^2(M),\nabla^{nc})$, $\beta \in Par(\Lambda^{k+1}(M),\nabla^{nc})$ and $g \in c$. Then it holds for the nc-Killing form $\left([\beta, \alpha]_{\mathcal{T}}\right)_+ \in \Omega_{nc,g}(M)$ that
\begin{align}
\frac{1}{2} \left([\beta, \alpha]_{\mathcal{T}}\right)_+ = L_{V_{\alpha}} \beta_+ - (k+1) \lambda_{\alpha} \cdot \beta_+ \in \Omega^k_{nc,g}(M). \label{equ}
\end{align}
Here, $L$ denotes the Lie derivative, $V_{\alpha}$ is the conformal vector field canonically associated to $\alpha$, and $\lambda_{\alpha} \in C^{\infty}(M)$ is defined via $L_{V_{\alpha}}g=2 \lambda_{\alpha} \cdot g$. In particular, the right hand side of (\ref{equ}) is again a nc-Killing $k$-form.
\end{Proposition}

\begin{bemerkung}
Proposition \ref{na} yields a natural action which gives the space of nc-Killing $k$-forms the structure of a module for the Lie algebra of normal conformal vector fields. In this context, we remark that it has already been shown in \cite{sem} that for a \textit{conformal} vector field $V$ and \textit{conformal} Killing $k-$form $\beta_+$, the form $L_{V} \beta_+ - (k+1) \lambda_{V} \cdot \beta_+$ is again a \textit{conformal} Killing $k-$form.
\end{bemerkung}

\textit{Proof. }
Dualizing the first nc-Killing equation (cf. \ref{stu}) for $\alpha_+$ yields
\begin{align}
\nabla^g_X V_{\alpha} = (X \invneg \alpha_0)^{\sharp} + \alpha_{\mp} X. \label{k45}
\end{align}
We have that $(L_{V_{\alpha}}g)(X,Y) = g(\nabla^g_X V_{\alpha},Y) + g(\nabla^g_Y V_{\alpha},X) = 2 \lambda_{\alpha} g(X,Y)$. Inserting (\ref{k45}) shows that $\alpha_{\mp}=\lambda_{\alpha} \in C^{\infty}(M)$. We fix $x \in M$ and let $(s_1,...,s_n)$ be a local $g-$pseudo-orthonormal frame around $x$. Cartans formula for the Lie derivative $L$ yields that around $x$ we have
\begin{align*}
L_{V_{\alpha}}\beta_+ &= d\left(V_{\alpha} \invneg \beta \right) + V_{\alpha} \invneg d \beta_+ \\
&=\sum_{i=1}^n \epsilon_i s_i^{\flat} \wedge \underbrace{\nabla^g_{s_i} \left(V_{\alpha} \invneg \beta_+\right)}_{=\left(\nabla^g_{s_i} V_{\alpha}\right) \invneg \beta_+ + V_{\alpha} \invneg \nabla^g_{s_i} \beta_+} + V_{\alpha} \invneg d \beta_+ \\
&= \sum_{i=1}^n \epsilon_i \left(s_i^{\flat} \wedge \left( \left(\nabla^g_{s_i} V_{\alpha}\right) \invneg \beta_+ \right) - V_{\alpha} \invneg \left( s_i^{\flat} \wedge \nabla^g_{s_i} \beta_+ \right) + g(s_i,V_{\alpha}) \cdot \nabla^g_{s_i} \beta_+ \right) + V_{\alpha} \invneg d \beta_+ \\
&= \underbrace{\sum_{i=1}^n \epsilon_i s_i^{\flat} \wedge \left( \left(\nabla^g_{s_i} V_{\alpha}\right) \invneg \beta_+ \right)}_{\text{I}} + \underbrace{\nabla^g_{V_{\alpha}} \beta_+}_{\text{II}}
\end{align*}
Using the nc-Killing equations for $\alpha_+$ and $\beta_+$, we rewrite the two summands as follows:
\begin{align*}
\text{I} =  \underbrace{\sum_{i=1}^n \epsilon_i s_i^{\flat} \wedge \left((s_i \invneg \alpha_0)^{\sharp} \invneg \beta_+ \right)}_{\text{I}a} + \underbrace{\alpha_{\mp} \cdot \sum_{i=1}^n \epsilon_i s_i^{\flat} \wedge \left(s_i  \invneg \beta_+ \right)}_{\text{I}b} 
\end{align*}
Clearly, I$b = k \cdot \alpha_{\mp} \cdot \beta_+$. In order to express I$a$ nicely, we introduce functions $a_{ij}$ such that $(s_i \invneg \alpha_0 )^{\sharp} = \sum_j \epsilon_j a_{ij} \cdot s_j$. Clearly, $a_{ij}=-a_{ji}$ and $\alpha_0 = \sum_{i < j} \epsilon_i \epsilon_j a_{ij} s_i^{\flat} \wedge s_j^{\flat}$. Inserting this into I$a$ yields that
\begin{align*}
\text{I}a = \sum_{i <j} \epsilon_i \epsilon_j a_{ij} \cdot \left(s_i^{\flat} \wedge \left(s_j \invneg \beta_+ \right) - s_j^{\flat} \wedge \left(s_i \invneg \beta_+ \right)\right). 
\end{align*}
In order to simplify this expression, we proceed as follows: Let $s^{\flat}_{J}:=s_{j_1}^{\flat} \wedge...\wedge s_{j_{k+1}}^{\flat}$ for $1 \leq j_1 < ... < j_{k+1} \leq n$. We compute for $i<j$:
\begin{align*}
s^{\flat}_J \cdot \left( s_i^{\flat} \wedge s_j^{\flat} \right) &= \left(s^{\flat}_J \cdot s_{i} \right) \cdot s_j = (-1)^{k+1}\left(s_i^{\flat} \wedge s_{J}^{\flat} + s_i \invneg s_J^{\flat} \right) \cdot s_j \\
&=s^{\flat}_i \wedge s^{\flat}_j \wedge s^{\flat}_J + s^{\flat}_i \wedge (s_j \invneg s^{\flat}_J ) - s^{\flat}_j \wedge (s^{\flat}_i \invneg s^{\flat}_J ) + s_i \invneg s_j \invneg s^{\flat}_J.
\end{align*}
Similarly, one obtains
\begin{align*}
\left( s_i^{\flat} \wedge s_j^{\flat} \right) \cdot s_J^{\flat} = s^{\flat}_i \wedge s^{\flat}_j \wedge s^{\flat}_J - s^{\flat}_i \wedge (s^{\flat}_j \invneg s^{\flat}_J ) + s^{\flat}_j \invneg (s^{\flat}_i \wedge s^{\flat}_J ) + s_i \invneg s_j \invneg s^{\flat}_J.
\end{align*}
Consequently, $\frac{1}{2} \cdot \left( s^{\flat}_J \cdot \left( s_i^{\flat} \wedge s_j^{\flat} \right) - \left( s_i^{\flat} \wedge s_j^{\flat} \right) \cdot s_J^{\flat} \right) = s_i^{\flat} \wedge \left(s_j \invneg s^{\flat}_J \right) - s_j^{\flat} \wedge \left(s_i \invneg s^{\flat}_J \right)$, and multilinear extension immediately yields that
\[\text{I}a =\frac{1}{2}(\beta_+ \cdot \alpha_0 - \alpha_0 \cdot \beta_+). \]
Furthermore, the summand II can with the nc-Killing equation for $\beta_+$ be rewritten as 
\begin{align*}
\nabla^g_{V_{\alpha}}\beta_+ &= V_{\alpha} \invneg \beta_0 + \alpha \wedge \beta_{\mp} \\
&=\frac{1}{2} \cdot \left((-1)^{k+1} \beta_0 \cdot \alpha_+ - \alpha_+ \cdot \beta_0 \right) + \frac{1}{2} \cdot \left( (-1)^{k+1} \beta_{\mp} \cdot \alpha_+ + \alpha_+ \cdot \beta_{\mp} \right).
\end{align*}
Putting all these formulas together again yields that
\begin{align*}
L_{V_{\alpha}} \beta_+ - (k+1) \lambda_{\alpha} \beta_+ =& \frac{1}{2} \left((-1)^{k+1}\beta_0 \cdot \alpha_+ - \alpha_+ \cdot \beta_0 +(-1)^{k+1} \beta_{\mp} \cdot \alpha_+ + \alpha_+ \cdot \beta_{\mp} + \beta_+ \cdot \alpha_0 - \alpha_0 \cdot \beta_+ \right) \\
&+ k \cdot \alpha_{\mp} \beta_+ -(k+1) \cdot \alpha_{\mp} \beta_+.
\end{align*}
Comparing this expression to (\ref{77}) immediately yields the Proposition.
$\hfill \Box$\\

As a second application of Proposition \ref{cne} we consider the case of $g$ being an Einstein metric in the conformal class. 

\begin{Proposition}
If $\beta \in \Omega_{nc,g}^k(M)$ is a nc-Killing $k$-form wrt. an Einstein metric $g$ on $M$, then both $\beta_0 = (k+1) \cdot d\beta_+$ and $\beta_{\mp} = -(n-k+1) \cdot d^*\beta_+$ are nc-Killing forms for $g$ as well.
\end{Proposition}

\textit{Proof. }As elaborated in \cite{nc}, on an Einstein manifold $(M,g)$, the tractor 1-form $\alpha = \left( 1,0,0,- \frac{\text{scal}^g}{2(n-1)n} \right)$ is parallel. Inserting this expression for $\alpha$ into the formulas in Proposition \ref{cne} shows that $\frac{1}{k+1}d\beta_+ + \frac{1}{n-k+1}d^*\beta_+$ is a nc Killing form.
$\hfill \Box$

\begin{bemerkung}
The last statement has a well-known spinorial analogue: Consider a twistor spinor $\ph \in \Gamma(S^g)$ on an Einstein manifold. As in this case $\nabla_X D^g \ph = \frac{n}{2} K^g(X) \cdot \ph = X \cdot \left(\frac{\text{scal}^g}{4n(n-1)}\cdot \ph \right)$, the spinor $D^g \ph$ turns out to be a twistor spinor on $(M,g)$ as well.
\end{bemerkung}

We compute the expression of the even-odd bracket wrt. a metric in the conformal class: 

\begin{Proposition} \label{prr}
Let $\alpha \in \Omega^{k+1}_{\mathcal{T}}(M)$ be a parallel tractor $(k+1)$-form, $\psi \in \Gamma(\mathcal{S}(M))$ a parallel spin tractor. For given $g \in c$ let $\Phi^g_{\Lambda}(\alpha) = (\alpha_+,\alpha_0,\alpha_{\mp},\alpha_-)$ with $\alpha_+ \in \Omega^k_{nc,g}(M)$ and $\widetilde{\Phi}^g(\psi) = (\ph, -\frac{1}{n}D^g\ph)$ with $\ph \in \text{ker }P^g$. Then the twistor spinor corresponding to the parallel spin tractor $[\alpha,\psi]=\alpha \cdot \psi \in \g_1$ via $g$ is given by 
\begin{align*}
\alpha_+ \circ \ph:=\widetilde{\Phi}^g ( proj_+^g \left(\alpha \cdot \psi \right)) &= \frac{2}{n} \alpha_+ \cdot D^g \ph + (-1)^{k+1} \alpha_{\mp} \cdot \ph + (-1)^{k+1} \alpha_{0} \cdot \ph \\
& = \frac{2}{n} \alpha_+ \cdot D^g \ph + \frac{(-1)^{k}}{n-k+1} d^*\alpha_{+} \cdot \ph + \frac{(-1)^{k+1}}{k+1} d\alpha_{+} \cdot \ph \in \text{ker }P^g.
\end{align*} 
\end{Proposition}

\textit{Proof. }
For given $x \in M$ we consider the reductions $\sigma^g:\mathcal{P}^g \rightarrow \overline{\mathcal{P}^1}$ and $\widetilde{\sigma}^g:\mathcal{Q}^g \rightarrow \overline{\mathcal{Q}^1}$ as introduced in chapter 2 with $\sigma^g \circ f^g = \overline{f}^1 \circ \widetilde{\sigma}^g$, and on some open neighbourhood $U$ of $x$ in $M$ we have 
\begin{align*}
\psi &= [\widetilde{\sigma}^g(\widetilde{u}), e_- \cdot w + e_+ \cdot w], \\
\alpha &= [\sigma^g(u),\underbrace{e_+^{\flat} \wedge \widetilde{\alpha}_+ + e_- \wedge e_+ \wedge \widetilde{\alpha}_{\mp} + \widetilde{\alpha}_0 + e_-^{\flat} \wedge \widetilde{\alpha}_-}_{ =: \widetilde{\alpha}}]
\end{align*}
for sections $\widetilde{u}:U \rightarrow \mathcal{Q}^g$, $u=f^g(\widetilde{u}):U \rightarrow \Pe^g$ and smooth functions $w: U \rightarrow \Delta_{p+1,q+1}$, $\widetilde{\alpha}_+,\widetilde{\alpha}_- : U \rightarrow \Lambda^k_{p,q}$, $\widetilde{\alpha}_0 : U \rightarrow \Lambda^{k+1}_{p,q}$ and $\widetilde{\alpha}_{\mp}: U \rightarrow \Lambda^{k-1}_{p,q}$. It follows by definition that on $U$
\begin{align*}
{\alpha} \cdot \psi = \left[\widetilde{\sigma}^g(\widetilde{u}), {\widetilde{\alpha}} \cdot (e_- \cdot w + e_+ \cdot w) \right].
\end{align*}
Consequently, we get for the corresponding twistor spinor wrt. $g$ that
\begin{align}
\widetilde{\Phi}^g ( {proj}_+^g \left({\alpha} \cdot \psi \right)) = \left[\widetilde{u}, \chi \left(e_- \cdot {proj}_{Ann(e_+)}\left({\widetilde{\alpha}} \cdot (e_- \cdot w + e_+ \cdot w)\right)\right) \right] \label{me}
\end{align}
Here, we identify the $Spin(p,q)-$modules $\Delta_{p,q}^{\C} \cong Ann(e_-)$ (cf. (\ref{fs})) by means of some fixed isomorphism $\chi$. One thus has to compute ${\widetilde{\alpha}} \cdot (e_- \cdot w + e_+ \cdot w)$. With the formulas for the action of $\Lambda_{p+1,q+1}^*$ on $\Delta_{p+1,q+1}$, it is straightforward to calculate that this product is given by
\begin{align*}
{\widetilde{\alpha}} \cdot (e_- \cdot w + e_+ \cdot w) =&(-1)^k \widetilde{\alpha}_+ \cdot e_+ \cdot e_- \cdot w + (-1)^k \widetilde{\alpha}_- \cdot e_- \cdot e_+ \cdot w + \widetilde{\alpha}_0 \cdot (e_- \cdot w + e_+ \cdot w) \\
&+ \widetilde{\alpha}_{\mp} \cdot (e_+ \cdot w - e_- \cdot w) \\
=& (e_- + e_+) \cdot ((-1)^k e_+ \cdot \widetilde{\alpha}_- \cdot w + (-1)^k e_- \cdot \widetilde{\alpha}_+ \cdot w + (-1)^{k+1} \widetilde{\alpha}_0 \cdot w \\
&+ (-1)^{k+1} \widetilde{\alpha}_{\mp} \cdot (e_+\cdot e_- \cdot w + w) ) \\
=:& (e_- + e_+) \cdot \widetilde{w}.
\end{align*}
Thus, one has by definition
\begin{align*}
\chi \left(e_- \cdot {proj}_{Ann(e_+)}\left({\widetilde{\alpha}} \cdot (e_- \cdot w + e_+ \cdot w)\right) \right) =& \chi \left(e_- \cdot e_+ \cdot \widetilde{w} \right) \\
=&-2 \widetilde{\alpha}_+ \cdot \chi(e_- \cdot w ) + (-1)^{k+1} \cdot \widetilde{\alpha}_0 \cdot \chi(e_- \cdot e_+ \cdot w )\\ 
&+ (-1)^{k+1} \cdot  \widetilde{\alpha}_{\mp} \cdot \chi(e_- \cdot e_+ \cdot w).
\end{align*}
Inserting this into (\ref{me}) yields that
\begin{align*}
\widetilde{\Phi}^g ( {proj}_+^g \left({\alpha} \cdot \psi \right)) =& -2\cdot[u,\widetilde{\alpha}_+] \cdot \underbrace{\left[\widetilde{u},\chi(e_- \cdot w)\right]}_{=\widetilde{\Phi}^g ({proj}_-^g \psi)}+ (-1)^{k+1} [u,\widetilde{\alpha}_{\mp}] \cdot \underbrace{[\widetilde{u},\chi(e_- \cdot e_+ \cdot w)]}_{=\widetilde{\Phi}^g ({proj}_+^g \psi)} \\
&+ (-1)^{k+1} [u,\widetilde{\alpha}_{0}] \cdot \underbrace{[\widetilde{u},\chi(e_- \cdot e_+ \cdot w)]}_{=\widetilde{\Phi}^g ({proj}_+^g \psi)} \\
=& \frac{2}{n} {\alpha}_+ \cdot D^g \ph + (-1)^{k+1}{\alpha}_{\mp} \cdot \ph + (-1)^{k+1} {\alpha}_{0} \cdot \ph.
\end{align*}
$\hfill \Box$

\begin{bemerkung}
In particular, Proposition \ref{prr} describes a principle of constructing new twistor spinors from a given twistor spinor and a nc-Killing form in an arbitrary pseudo-Riemannian setting. One can also show independently and more directly, i.e. without using tractor calculus, that for a given nc-Killing form $\alpha_+ \in \Omega^k_{nc,g}(M)$ and $\ph \in$ ker $P^g$, the spinor 
\begin{align}
\alpha_+ \circ \ph := \frac{2}{n} \alpha_+ \cdot D^g \ph + \frac{(-1)^{k}}{n-k+1} d^*\alpha_{+} \cdot \ph + \frac{(-1)^{k+1}}{k+1} d\alpha_{+} \cdot \ph \in \Gamma(S^g) \label{fo2}
\end{align}
 is again a twistor spinor on $(M,g)$.
To this end, we compute $\nabla^{S^g}_X (\alpha_+ \circ \ph)$ for $X \in \mathfrak{X}(M)$ using the nc-Killing formulas (cf. (\ref{stu})):
\begin{align*}
\nabla^{S^g}_X \left( \alpha_+ \cdot D^g \ph \right) & = \left( \nabla^g_X \alpha_+ \right) \cdot D^g \ph + \alpha_+ \cdot \nabla^{S^g}_X D^g \ph \\
 &= (X \invneg \alpha_0) \cdot D^g \ph + (X^{\flat} \wedge \alpha_{\mp}) \cdot D^g \ph + \alpha_+ \cdot \left(\frac{n}{2} \cdot K^g(X) \cdot \ph \right), \\
\nabla^{S^g}_X \left( \alpha_0 \cdot \ph \right) & = \left( \nabla^g_X \alpha_0 \right) \cdot \ph + \alpha_0 \cdot \nabla^{S^g}_X \ph \\
&= (K^g(X) \wedge \alpha_+) \cdot \ph - (X^{\flat} \wedge \alpha_-) \cdot \ph - \frac{1}{n} \cdot \alpha_0 \cdot X \cdot D^g \ph, \\
\nabla^{S^g}_X \left( \alpha_{\mp} \cdot \ph \right) & = \left( \nabla^g_X \alpha_{\mp} \right) \cdot \ph + \alpha_{\mp} \cdot \nabla^{S^g}_X \ph \\
&= (K^g(X) \invneg \alpha_+) \cdot \ph + (X  \invneg \alpha_-) \cdot \ph - \frac{1}{n} \cdot \alpha_{\mp} \cdot X \cdot D^g \ph. \\
\end{align*}
We deduce using the formulas (\ref{ext1}) that $\nabla^{S^g}_X (\alpha_+ \circ \ph) = X \cdot \xi$ for all $X \in \mathfrak{X}(M)$, where $\xi := \left(\frac{1}{n} \alpha_0 \cdot D^g \ph + \frac{1}{n} \alpha_{\mp} \cdot D^g \ph + (-1)^{k+1} \alpha_- \cdot \ph \right)$, showing that $\alpha_+ \circ \ph$ satisfies the twistor equation with $D^g(\alpha_+ \circ \ph)=-n \cdot \xi$.
\end{bemerkung}

Finally, we discuss the case of $\alpha_+$ being a nc-Killing 1-form and $V_{\alpha}$ the dual normal conformal vector field\footnote{The proof of the following statement is then also the postponed proof of Proposition \ref{sld}.}.

\begin{Proposition} \label{lsd}
In the setting of Proposition \ref{prr}, if $k=1$ we have
\begin{align}
\Phi^g \left( proj_+^g \left(\alpha \cdot \psi \right)\right) = -2 \cdot \underbrace{\left(\nabla_{V_{\alpha}} \ph + \frac{1}{4} \tau \left(\nabla V_{\alpha} \right) \cdot \ph  \right)}_{=:V_{\alpha} \circ \ph}, \label{sld2}
\end{align}
where $\tau \left(\nabla V_{\alpha} \right) = \sum_{j=1}^n \epsilon_j \left( \nabla_{s_j} V_{\alpha} \right) \cdot s_j + (n-2) \cdot \lambda_{\alpha}$ for any local $g-$pseudo-orthonormal frame $(s_1,...,s_n)$, and $L_{V_{\alpha}}g = 2 \lambda_{\alpha}g$. 
\end{Proposition}
\textit{Proof. }Wrt. $g$ it holds that $\Phi^g(\alpha) = (\alpha_+,\alpha_0,\alpha_{\mp},\alpha_-)$. As in the proof of Proposition \ref{na} it follows that $\alpha_{\mp} = \lambda_{\alpha}$. Let $(s_1,...,s_n)$ be $g-$orthonormal. The first nc-Killing equation for $\alpha_+$ yields that $\nabla^g_{s_j} V_{\alpha} = \left(s_j \invneg \alpha_0\right)^{\sharp} + \alpha_{\mp} \cdot s_j$. Right-multiplication by $s_j$ gives $\left( \nabla^g_{s_j} V_{\alpha} \right) \cdot s_j = - (s_j \wedge (s_j \invneg \alpha_0)) - \epsilon_j \cdot \alpha_{\mp}$. Summing over $j$ thus reveals that $\tau \left(\nabla V_{\alpha} \right) = -2 \cdot \alpha_0 - n \cdot \alpha_{\mp}$, and together with the twistor equation we conclude that the right-hand side of (\ref{sld2}) is given by
\begin{align*}
-2 \cdot \left( -\frac{1}{n} V_{\alpha} \cdot D^g \ph - \frac{1}{2} \alpha_0 \cdot \ph - \frac{1}{2} \alpha_{\mp} \cdot \ph \right).
\end{align*}
Comparing this to the result of Proposition \ref{prr} immediately yields (\ref{sld2}).
$\hfill \Box$

\begin{bemerkung}
The term $V_{\alpha} \circ \ph$ in (\ref{sld2}) has become standard in the literature as \textsf{spinorial Lie derivative} as introduced in \cite{kos,ha96,raj}. Thus, the metric description of the even odd bracket in Proposition \ref{prr} can be viewed as a \textsf{generalization of the spinorial Lie derivative} to higher order nc-Killing forms, and we see that the brackets in the tractor conformal superalgebra reproduce the spinorial Lie derivative when a metric is fixed. For the case $k=1$, \cite{ha96} shows that $X \circ \ph$ is a twistor spinor for every twistor spinor $\ph$ and every \textit{conformal} vector field $X$, i.e. $X$ need not to be normal conformal.
\end{bemerkung}

\begin{bemerkung}
As in the Lorentzian setting, it is also possible in arbitrary signatures to include all conformal Killing forms, i.e. not only nc-Killing forms, in the even part of the algebra in terms of distinguished tractors. However, the generalization of (\ref{star}) to arbitrary signatures, which can be found in \cite{cov2}, is technically very demanding.
\end{bemerkung}

\subsection*{Analogous construction for special Killing forms and Killing spinors}
We specialize the principle for constructing new nc-Killing forms out of existing ones using the $\circ-$operations from Proposition \ref{cne}. In this context, we make some more general definitions and remarks:

\begin{definition}
Let $(M^{p,q},g)$ be a pseudo-Riemannian manifold of constant scalar curvature $\text{scal}^g$. A $k-$form $\alpha \in \Omega^k(M)$ is called a special Killing k-form to the Killing constant $-\frac{(k+1)\text{scal}^g}{n(n-1)}$ if
\begin{equation} \label{fgt}
\begin{aligned}
\nabla_X^g \alpha &= \frac{1}{k+1} X \invneg d \alpha, \\
\nabla^g_X d\alpha &= - \frac{(k+1) \text{scal}^g}{n(n-1)} \cdot X^{\flat} \wedge \alpha.
\end{aligned}
\end{equation}
We let $\Omega^k_{sk,g}(M)$ denote the space of all special Killing $k-$forms on $(M,g)$.
\end{definition}

Examples and classification results for special Killing forms are discussed in \cite{sem}. For instance, the dual of every Killing vector field defining a Sasakian structure and the Dirac currents of real Killing spinors on Riemannian manifolds are special Killing 1-forms. Note that every special Killing form is conformal and coclosed, i.e. $d^* \alpha = 0$.\\
\newline
Let us from now on assume that $\text{scal}^g \neq 0$.  Under this assumption, spaces carrying special conformal Killing forms can be classified using an analogue of B\"ars cone construction for Killing spinors, see \cite{baer}, for differential forms. More precisely, consider the cone $C(M)=\R^+ \times M$ with cone metric $\widehat{g}_b := bdt^2 +t^2g$, where $b \neq 0$ is a constant scaling, of signature $(p,q+1)$ or $(p+1,q)$. 

\begin{Proposition}{\cite{sem}} \label{pse}
Let $b=\frac{(n-1)n}{\text{scal}^g}$. Then special Killing $k-$forms to the Killing constant $-\frac{(k+1)\text{scal}^g}{n(n-1)}$ are in 1-to-1 correspondence to parallel $(k+1)$-forms on the cone $(C(M), \widehat{g}_b)$, given by
\begin{align}
\Omega_{sk,g}^k(M) \ni \alpha \leftrightarrow \widehat{\alpha}:= t^k dt \wedge \alpha + \frac{\text{sgn}(b) t^{k+1}}{k+1} d \alpha \in \Omega^{k+1}(C(M)) \label{sug}
\end{align}
\end{Proposition}
Using this, one classifies compact, simply-connected Riemannian manifolds carrying special Killing forms, see \cite{sem}. We come back to this list in the last section of this thesis.

\begin{bemerkung}
One can now derive analogous formulas to (\ref{fo1}),(\ref{fo2}) for special Killing forms and Killing spinors on pseudo-Riemannian manifolds using the cone construction, i.e. proceed as follows:
\begin{enumerate}
\item We let $\alpha \in \Omega^k(M), \beta \in \Omega^l(M)$ be special Killing forms to the same Killing constant and $\ph \in \Gamma(S^g)$ a Killing spinor on $(M,g)$.
\item Using B\"ars construction and Proposition \ref{pse}, we view these objects as parallel tensors $\widehat{\alpha}, \widehat{\beta} \in \Omega^{k+1}(C(M))$, $\widehat{\beta} \in \Omega^{l+1}(C(M))$ and $\widehat{\ph} \in \Gamma(C(M),S^{\widehat{g}_b})$.
\item We compute $\widehat{\alpha} \cdot \widehat{\beta}$ (with (\ref{prodo}) applied pointwise) and $\widehat{\alpha} \cdot \widehat{\ph}$ which again turn out to be parallel forms resp. spinors on the cone.
\item Via (\ref{sug}), one expresses these products as special Killing forms resp. Killing spinors  on the base $(M,g)$ using the original data $\alpha, \beta, \ph$ and $d \alpha$, $d\beta$ only. Let us call these objects $\alpha \circ \beta \in \Omega^*_{sk,g}(M)$ and $\alpha \circ \ph \in \mathcal{K}(M)$.
\end{enumerate}
Carrying these steps out is straightforward. One obtains the same formulas (\ref{fo1}) and (\ref{fo2}), which of course simplify since $d^*\alpha = 0, D^g \ph = -\lambda \cdot n \cdot \ph$ for some $\lambda \in i\R \cup \R$ with $\ph \in \mathcal{K}_{\lambda}(M)$.  In other words, one obtains a map
\begin{equation} \label{fo21}
\begin{aligned}
\circ: \Omega^k_{sk,g}(M) \times \Omega^l_{sk,g}(M) & \rightarrow \Omega^*_{sk,g}(M), \\
(\alpha, \beta)  & \mapsto \alpha \circ \beta = \frac{1}{l+1}\cdot \alpha \cdot d\beta + (-1)^{k+1} \frac{1}{k+1} d\alpha \cdot \beta,
\end{aligned}
\end{equation}
and an action of special Killing forms on Killing spinors, given by
\begin{equation} \label{fo22}
\begin{aligned}
\circ: \Omega^k_{sk,g}(M) \times \mathcal{K}_{\lambda}(M) &\rightarrow \mathcal{K}_{\lambda} \oplus \mathcal{K}_{- \lambda}(M), \\
(\alpha, \ph) & \mapsto \alpha \circ \ph = -2 \cdot \alpha \cdot \ph + \frac{(-1)^{k+1}}{k+1} d\alpha \cdot \ph.
\end{aligned}
\end{equation}
In particular,  (\ref{fo21}) allows one to construct new special Killing forms out of existing special Killing forms.
\end{bemerkung}

However, for pseudo-Riemannian Einstein spaces which are not Ricci-flat, special Killing forms are more directly related to normal conformal Killing forms and there is an equivalent way of deriving (\ref{fo21}) and (\ref{fo22}):\\
\cite{leihabil,lst} shows that for every pseudo-Riemannian Einstein space $(M,g)$, the conformal holonomy coincides with the holonomy of an ambient space which is the cone trivially extended by a parallel direction, i.e. $Hol(M,[g])=Hol(C(M),\widehat{g}_b)$. Using this, it is easy to deduce that there is a natural and bijective correspondence between parallel tractor forms on $M$, i.e. normal conformal Killing forms for $(M,g)$, and parallel forms on the cone, i.e. special Killing forms for $(M,g)$. More precisely, one shows:

\begin{Proposition}[\cite{leihabil}] \label{sumsibum}
On a pseudo-Riemannian Einstein space of nonvanishing scalar curvature, every nc-Killing form is the sum of a special Killing form and a closed Killing form.
\end{Proposition}
In particular, the coclosed nc-Killing forms on Einstein spaces are precisely the special Killing forms. This also follows from a direct inspection of the nc-Killing equations. Note that the well-known spinorial analogue of Proposition \ref{sumsibum} is the fact that on an Einstein space every twistor spinor decomposes into the sum of two Killing spinors. Thus, for Einstein spaces one obtains the maps (\ref{fo21}) and (\ref{fo22}) by restriction of (\ref{map}) and (\ref{fo2}) to special Killing forms and Killing spinors.

\section{The possible dimensions of the space of twistor spinors} \label{podi}
We have already discussed for Lorentzian signatures, in how far algebraic structures of the tractor conformal superalgebra $\g = \g_0 \oplus \g_1$, in particular, whether it is a \textit{Lie} superalgebra, are related to (local) geometric structures in the conformal class. It is natural to investigate this question further in arbitrary signatures, and we ask ourselves how possible dimensions of the odd \textit{supersymmetric} part $\g_1$ are related to underlying geometries. Main ingredient is the following algebraic Lemma:
\begin{Lemma} \label{est}
For integers $r$ and $s$ consider the bilinear map 
\begin{align*}
V: \Delta_{r,s}^{\R} \otimes \Delta_{r,s}^{\R} \rightarrow \R^{r,s}\text{, }(\psi_1,\psi_2) \mapsto V_{\psi_1,\psi_2}
\end{align*}
mapping a pair of spinors to the associated vector. Let $S_0 \subset \Delta^{\R}_{r,s}$ be a linear subspace and set $V_{S_0}:=V_{|S_0 \otimes S_0 }$. We have:
\begin{enumerate}
\item If dim $S_0 > \frac{3}{4}\cdot\text{dim } \Delta_{r,s}^{\R}$, then $V_{S_0}$ is surjective.
\item If dim $S_0 > \frac{1}{2}\cdot\text{dim } \Delta_{r,s}^{\R}$, then $V_{S_0}$ is not the zero map.
\end{enumerate}
\end{Lemma}

\textit{Proof. }The first part is proved in  \cite{cortes}. For the second part, assume that $V_{S_0}(\psi_1,\psi_2)=0$ for all $\psi_1,\psi_2 \in S_0$. By definition, this is equivalent to $\langle v \cdot \psi_1 , \psi_2 \rangle_{\Delta_{r,s}^{\R}} = 0$ for all $\psi_1,\psi_2 \in S_0$ and $v \in \R^{r,s}$, i.e. $\forall v \in \R^{r,s}: cl(v): S_0 \rightarrow S_0^{\bot}$. As dim $S_0 > \frac{1}{2} \cdot \text{dim }\Delta_{r,s}^{\R}$, it follows that dim $S_0^{\bot} < \frac{1}{2} \cdot \text{dim }\Delta_{r,s}^{\R}$. Thus the map $cl(v)$ has a kernel for every $v \in \R^{r,s}$, i.e. there is $\psi_v \in \Delta^{\R}_{r,s} \backslash \{0 \}$ with $v \cdot \psi_v = 0$. This implies that $\langle v,v \rangle_{r,s} = 0$ for every $v \in \R^{r,s}$. 
$\hfill \Box$

\begin{bemerkung}
The second statement in Lemma \ref{est} cannot be improved in general. Namely, taking $r=s=2$ and $S_0 := \Delta_{2,2}^{\R,\pm} \subset \Delta_{2,2}^{\R}$ provides an example for dim $S_0 = \frac{1}{2} \cdot \text{dim } \Delta_{r,s}^{\R}$ and $V_{S_0} = 0$.
\end{bemerkung}

Applications of Lemma \ref{est} have already been studied in the literature::

\begin{Proposition}{\cite{cortes}} \label{copo}
Let $(M^{p,q},g)$ be a pseudo-Riemannian spin manifold of dimension $n$ with real spinor bundle $S^g = S^g_{\R}(M)$ of rank $N$. 
\begin{enumerate}
\item If $(M,g)$ admits $k > \frac{3}{4}N$ twistor spinors which are linearly independent at $x \in M$, then $(M,g)$ admits $n$ conformal vector fields, which are linearly independent at $x \in M$.
\item If $(M,g)$ admits $k > \frac{3}{4}N$ parallel spinors, then $(M,g)$ is flat.
\end{enumerate}
\end{Proposition}

We now apply Lemma \ref{est} in the tractor setting yielding a conformal analogue of the second part of Proposition \ref{copo}. Let $(M^{p,q},c)$ be a conformal spin structure with real spin tractor bundle $\mathcal{S}(M)$ and space of \textit{real} twistor spinors $\g_1$. Let $N_c:= 2 \cdot \text{dim } \Delta_{p,q}^{\R}$ denote the rank of $\mathcal{S}(M)$, which is the maximal number of linearly independent real twistor spinors on $(M,c)$.

\begin{Proposition} \label{stuv}
In the above notation, we have:
\begin{enumerate}
\item If dim $\g_1 > \frac{3}{4} \cdot N_c$, then $(M,c)$ is conformally flat.
\item If dim $\g_1 > \frac{1}{2} \cdot N_c$, then there exists an Einstein metric in $c$ (at least on an open and dense subset).
\end{enumerate}
\end{Proposition}

\textit{Proof. }We apply Lemma \ref{est} to the case that $r=p+1$, $s=q+1$ and $S_0\subset \Delta_{p+1,q+1}^{\R}$ being the subspace of $Hol(M,c)$-invariant spinors (for some fixed base points) which as we know correspond to twistor spinors. Surjectivity of $V$ yields a basis of $\R^{p+1,q+1}$ which is $Hol(M,c)$-invariant. This proves the first part.\\
For the second part, it follows analogously by nontriviality of $V_{S_0}$ that there exists at least one nontrivial holonomy-invariant vector. By \cite{lst} this yields an Einstein scale in the conformal class (on an open, dense subset).
$\hfill \Box$

\begin{bemerkung}
As a simply-connected, conformally flat manifold always admits the maximal number of twistor spinors, the previous Proposition implies that either dim $\g_1 \leq \frac{3}{4} \cdot N_c$ or dim $g_1 = N_c$ is maximal, i.e. the dimension of $\g_1$ cannot be arbitrary for the simply-connected case.
\end{bemerkung}
In the second case of Proposition \ref{stuv} one can say more: To this end, let $(M^n,c=[g])$ be a  simply-connected pseudo-Riemannian conformal spin manifold where $g$ is a Ricci-flat metric. Let $k$ denote the number of linearly independent parallel vector fields on $(M,g)$. By \cite{lst} we have for $x \in M$ that in the Ricci-flat case
\begin{align*}
\mathfrak{hol}_x(M,[g]) = \mathfrak{hol}_x(M,g) \ltimes {\R}^{n-k} = \left\{ \begin{pmatrix} 0 & v^{\flat} & 0 \\ 0 & A & -v \\ 0 & 0 & 0 \end{pmatrix} \mid A \in  \mathfrak{hol}_x(M,g), v \in \R^{n-k} \right\},
\end{align*}
where the matrix is written wrt. the basis $(s_+,s_1,...,s_n,s_-)$ of $\mathcal{T}_x(M)$ for some pseudo-orthonormal basis $(s_1,...,s_n)$ of $T_xM$. Assume now that $k < n$, i.e. $(M,g)$ is Ricci-flat but non flat, and let $\psi \in \g_1$ be a parallel spin tractor on $(M,[g])$ with twistor spinor $\ph$. It follows by the holonomy-principle that
\begin{align}
\lambda_*^{-1} \left(  \begin{pmatrix} 0 & v^{\flat} & 0 \\ 0 & 0 & -v \\ 0 & 0 & 0 \end{pmatrix} \right) \cdot \psi(x) = 0, \label{disco}
\end{align}
for all $v \in \R^{n-k} \subset \R^n$, i.e. $s_+ \cdot v \cdot \psi(x) = 0$ (cf. \cite{mat} for formulas for $\lambda_*^{-1}$ in this situation). Choosing $v$ to be non lightlike yields that $s_+ \cdot \psi(x) = 0$ for all $x \in M$ which is equivalent to $D^g \ph = -n \cdot \widetilde{\Phi}^g({proj}^g_- \psi) = 0$. Thus, $\ph$ is a parallel spinor on $(M,g)$, ker $\psi \neq \{0 \}$, and we have proved:

\begin{Proposition} \label{let}
Let $(M,g)$ be a simply-connected Ricci-flat spin manifold. Then either every twistor spinor on $(M,g)$ is parallel or $(M,g)$ is flat. In particular, if for a conformal structure $(M,c)$ there is a Ricci-flat metric in the conformal class and dim $\mathfrak{g}_1$ is not maximal, then $\mathfrak{g}$ is a Lie superalgebra.
\end{Proposition}
We now come back to the second case of Proposition \ref{stu}: It follows now directly from Proposition \ref{let} that in case of dim $\g_1 > \frac{1}{2} \cdot {\left( 2 \cdot \text{dim} \Delta_{p,q}^{\R}\right)}$ there exists an Einstein metric in $c$ with nonzero scalar curvature or the conformal structure is conformally flat, provided that $M$ is simply-connected.

\begin{beispiel}
We consider a special class of conformally Ricci-flat Lorentzian metrics admitting twistor spinors, namely plane waves $(M,h)$ which are equivalently characterized by the existence of local coordinates $(x,y_1,...,y_n,z)$ such that 
\begin{align*}
h = 2dx dz + \left(\sum_{i,j=1}^n a_{ij}y_i y_j \right) dz^2 + \sum_{i=1}^n dy_i^2,
\end{align*}
where the $a_{ij}$ are functions only of $z$. It is $Ric^h = \sum_{i=1}^n a_{ii} dz^2$ and the isotropic vector field $\frac{\partial}{\partial x}$ is parallel. Let us assume that $(M,h)$ is indecomposable. Then it is known from \cite{lcc} that for $x$ in $M$
\begin{align*}
\mathfrak{hol}_x(M,[h]) = \R^{2n+1} = \left\{ \begin{pmatrix} 0 & 0 & u^T & c & 0 \\ 0 & 0 & v^T & 0 & -c \\ 0 & 0 & 0 & -v & -u \\ 0 & 0 & 0 & 0 & 0 \\ 0 & 0 & 0 & 0 & 0 \end{pmatrix} \mid u,v \in \R^n, c \in \R \right\}
\end{align*}
This explicit description makes it straightforward to calculate all spinors annihilated by $\lambda_*^{-1}\left(\mathfrak{hol}_x(M,[h])\right)$, yielding that dim ker $P^g = \frac{1}{2}\cdot $ dim $\Delta^{\R}_{1,n+1} = \frac{1}{4} \cdot N_c$.
\end{beispiel}


\textbf{Acknowledgement} The author gladly acknowledges support from the DFG (SFB 647 - Space Time Matter at Humboldt University Berlin) and the DAAD (Deutscher Akademischer Austauschdienst / German Academic Exchange Service). Furthermore, it is a pleasure to thank Jose Figueroa O Farrill for various discussions about mathematical phsycis.

\small
\bibliographystyle{plain}
\bibliography{literatur}

\medskip

\medskip

\textsc{Andree Lischewski\\
Humboldt-Universit\"at zu Berlin, Institut f\"ur Mathematik\\
Rudower Chaussee 25, Room 1.310, D12489 Berlin.\\
E-Mail: }\texttt{lischews@mathematik.hu-berlin.de}

\end{document}